\documentclass[11pt,a4paper]{article}
\usepackage{amsfonts,amsmath,amssymb,accents,bm,amsthm}
\usepackage{verbatim}

\usepackage{hyperref}
\usepackage[normalem]{ulem}

\usepackage{tikz} 
\usetikzlibrary{math,calc,arrows,positioning,fit,petri,external}

\usepackage{graphicx}

\newcommand{\dis}{\displaystyle}
\newcommand{\txt}{\textstyle}

\newcommand{\noi}{\noindent}
\newcommand{\halmos}{\rule{1ex}{1.4ex}}

\newcommand{\QED}{\nopagebreak{\hspace*{\fill}$\halmos$\medskip}}
\newcommand{\med}{\medskip}
\newcommand{\quand}{\quad\mbox{and}\quad}

\newtheoremstyle{mythm}
  {}
  {}
  {\itshape}
  {}
  {\bfseries}
  {}
  {.5em}
  {#1 #2 \thmnote{(#3)}}

\theoremstyle{mythm}
\newtheorem{theorem}{Theorem}[section]
\newtheorem{proposition}[theorem]{Proposition}
\newtheorem{lemma}[theorem]{Lemma}
\newtheorem{conjecture}[theorem]{Conjecture}

\newcommand{\bt}{\begin{theorem}}
\newcommand{\et}{\end{theorem}}
\newcommand{\bl}{\begin{lemma}}
\newcommand{\el}{\end{lemma}}
\newcommand{\bp}{\begin{proposition}}
\newcommand{\ep}{\end{proposition}}
\newcommand{\bcon}{\begin{conjecture}}
\newcommand{\econ}{\end{conjecture}}

\newenvironment{Proof}[1][]{\noi\textbf{Proof #1}}{\QED}
\newcommand{\bpro}{\begin{Proof}}
\newcommand{\epro}{\end{Proof}}

\newcommand{\be}{\begin{equation}}
\newcommand{\ee}{\end{equation}}
\newcommand{\ba}{\begin{array}}
\newcommand{\ea}{\end{array}}
\newcommand{\bac}{\begin{array}{r@{\,}c@{\,}l}}
\newcommand{\bc}{\be\begin{array}{r@{\,}c@{\,}l}}
\newcommand{\ec}{\end{array}\ee}

\newcommand{\al}{\alpha}

\newcommand{\ga}{\gamma}
\newcommand{\Ga}{\Gamma}
\newcommand{\de}{\delta}
\newcommand{\De}{\Delta}
\newcommand{\eps}{\varepsilon}
\newcommand{\la}{\lambda}
\newcommand{\La}{\Lambda}
\newcommand{\sig}{\sigma}

\newcommand{\tet}{\theta}
\newcommand{\vtet}{\vartheta}
\newcommand{\om}{\omega}

\newcommand{\si}{\ensuremath{\sigma}}

\newcommand{\Ai}{{\cal A}}

\newcommand{\Fi}{{\cal F}}
\newcommand{\Gi}{{\cal G}}

\newcommand{\Si}{{\cal S}}
\newcommand{\Ti}{{\cal T}}

\newcommand{\Xc}{{\cal X}}

\newcommand{\E}{{\mathbb E}}

\newcommand{\N}{{\mathbb N}}

\renewcommand{\P}{{\mathbb P}}

\newcommand{\R}{{\mathbb R}}

\newcommand{\Z}{{\mathbb Z}}

\newcommand{\Xb}{{\mathbf X}}
\newcommand{\Yb}{{\mathbf Y}}

\newcommand{\rb}{{\mathbf r}}

\newcommand{\volgt}{\ensuremath{\Rightarrow}}
\newcommand{\up}{\uparrow}
\newcommand{\down}{\downarrow}
\newcommand{\sub}{\subset}
\newcommand{\beh}{\backslash}

\newcommand{\Asto}[1]{\underset{{#1}\to\infty}{\Longrightarrow}}

\newcommand{\To}{\Rightarrow}

\newcommand{\ti}{\tilde}

\newcommand{\dgg}{\dagger}
\newcommand{\ov}{\overline}
\newcommand{\un}{\underline}

\newcommand{\ffrac}[2]{{\textstyle\frac{{#1}}{{#2}}}}
\newcommand{\dif}[1]{\ffrac{\partial}{\partial{#1}}}

\newcommand{\cn}{\colon}
\newcommand{\di}{\mathrm{d}}
\newcommand{\half}{{[0,\infty)}}
\newcommand{\expo}{\mbox{\large\it e}}
\newcommand{\ex}[1]{\expo^{\,\textstyle{#1}}}
\newcommand{\ha}{\ffrac{1}{2}}

\newcommand{\cov}{{\rm Cov}}

\DeclareMathOperator*{\bplim}{\mbox{\rm bp-lim}}
\DeclareMathOperator*{\mbplim}{\mbox{\rm -bp-lim}}
\DeclareMathOperator*{\pbplim}{\mbox{\rm pbp-lim}}

\newcommand{\Lt}{L}

\begin{document}

\makeatletter\@addtoreset{equation}{section}
\makeatother\def\theequation{\thesection.\arabic{equation}}

\renewcommand{\labelenumi}{{\rm (\roman{enumi})}}
\renewcommand{\theenumi}{\roman{enumi}}

\title{\vspace*{-1cm}Equivalence of finite and infinite survival and sharpness of the phase transition for additive particle systems}
\author{Concetta Campailla\footnote{Universität Innsbruck, Institut für Mathematik, Technikerstraße 13,
A-6020 Innsbruck, Austria. 
maria.campailla@uibk.ac.at}\ \ and Jan~M.~Swart\footnote{The Czech Academy of Sciences, Institute of Information Theory and Automation, Pod vod\'arenskou v\v{e}\v{z}\'i~4, 18200 Praha 8, Czech Republic. swart@utia.cas.cz}} 

\date{\today}

\maketitle

\begin{abstract}\noi
Krone's two-stage contact process is an interacting particle system in which each site on a grid can be in one of three possible states: 0 (empty), 1 (young individual), and 2 (adult individual). Only adults can multiply, and their offspring is initially young, before growing up after an exponential time. There are two death rates, one for adult trees, and a higher one for young trees. The process is additive and has an additive dual, but unlike the contact process it is not self-dual. Answering a question of Krone, for the process on the $d$-dimensional integer grid, Foxall proved that survival of finite populations is equivalent to nontriviality of the upper invariant law. Using a completely different method, we generalise this result to processes on countable groups and along the way also establish sharpness of the phase transition for a general class of additive interacting particle systems.
\end{abstract}

\noi
%
%

\noi

{\setlength{\parskip}{-2pt}\tableofcontents}

\newpage

\section{Introduction}

Krone's two-stage contact process \cite{Kro99} is a Markov process $(X_t)_{t\geq 0}$ whose state space is the space of all functions $x$ from a countable grid $\La$ into the set $\{0,1,2\}$. Here $X_t(i)=0,1,2$ means that at time $t\geq 0$, the site $i\in\La$ is empty, occupied by a young tree, and occupied by an adult tree, respectively. An adult tree at $i$ gives birth with rate $\la(i,j)$ to a young tree at $j$, provided $j$ was previously empty. Because of the similarity to the contact process, we call $\la(i,j)$ the \emph{infection rate} from $i$ to $j$. Young trees mature at rate $\ga$, adult trees die at rate $\de_2$, and young trees die at rate $\de_1+\de_2$. The fact that young trees die at a higher rate than adult trees represents the biological reality that in many species, immature individuals have a higher death rate than adults.

Like the contact process, the two-stage contact process is additive and hence has an additive dual. Unlike the contact process, however, it is not self-dual: Krone identified its dual as the \emph{on-off contact process}, which differs from the two-stage contact process if $\de_1>0$. This means that the usual proof that for the contact process, finite and infinite survival are equivalent, breaks down for the two-stage contact process. This led Krone to ask \cite[page 349]{Kro99}:
\begin{quote}
\emph{Are finite survival and infinite survival equivalent for the two-stage
contact process with $\de>0$ or are there two distinct critical values?}
\end{quote}
For the process on $\La=\Z^d$ with infection rates that are finite range and translation invariant, Foxall \cite[Thm~1.2]{Fox15} gave a positive answer to Krone's question. His proof depends on the block construction of \cite{BG90,BG94} and is not easily adapted to other grids.

In this paper, we generalise Foxall's result by dropping the finite range condition and allowing the grid $\La$ to be any countable group, assuming only that the infection rates are invariant under the left action of the group on itself (see Theorem \ref{T:fininf}). Our method of proof is completely different from the one in \cite{Fox15}. We do not rely on \cite{BG90,BG94} but instead rely on the fact that for an additive system and its dual, in the subcritical regime, the expected number of infected sites decays at a rate that is the same for the forward and dual process, see Lemmas \ref{L:exprate} and \ref{L:duexp} below. This allows us to answer Krone's question by proving sharpness of the phase transition for a class of additive particle systems that includes the two-stage contact process and its dual. We do so by adapting the method of \cite{Swa18} to a general class of additive particle systems, culminating in Theorem~\ref{T:sharp2} below, which is the second main result of our paper.

\section{Statement of the results}

\subsection{Definition of the model}\label{S:2stage}

Let $\La$ be a countable group with group operation $(i,j)\mapsto ij$, inverse operation $i\mapsto i^{-1}$, and unit element $0$, also referred to as the origin. Let $\Si=\Si(\La):=\{0,1,2\}^\La$ be the space of all functions $x\cn\La\to\{0,1,2\}$. We set $\Si_{\rm fin}=\Si_{\rm fin}(\La):=\{x\in\Si:\sum_{i\in\La}x(i)<\infty\}$. For each $i\in\La$, we define maps from $\Si$ into itself by
\be\ba{l}\label{mapi}
\dis m^{1\to 2}_i(x)(k):=\left\{\ba{@{}l@{\hspace{7pt}}l}
2&\mbox{if }k=i,\ x(i)=1,\\
x(k)&\mbox{otherwise,}\ea\right.
\quad
m^{1\to 0}_i(x)(k):=\left\{\ba{@{}l@{\hspace{5pt}}l}
0&\mbox{if }k=i,\ x(i)=1,\\
x(k)&\mbox{otherwise,}\ea\right.\\[15pt]
\dis m^{2\to 1}_i(x)(k):=\left\{\ba{@{}l@{\hspace{5pt}}l}
1&\mbox{if }k=i,\ x(i)=2,\\
x(k)&\mbox{otherwise,}\ea\right.
\quad
m^{\to 0}_i(x)(k):=\left\{\ba{@{}l@{\hspace{5pt}}l}
0&\mbox{if }k=i,\\
x(k)&\mbox{otherwise.}\ea\right.
\ec
For each $i,j\in\La$, we moreover define
\be\label{mapij}
m^\ast_{i,j}(x)(k):=\left\{\ba{@{\,}l@{\hspace{10pt}}l}
1&\mbox{if }k=j,\ x(i)=2,\ x(j)=0,\\
x(k)&\mbox{otherwise.}\ea\right.
\ee
We fix $\ga>0$ and $\de_1,\de'_1,\de_2\geq 0$ as well as a function $\la\cn\La\times\La\to\half$ such that $\la(i,i)=0$ $(i\in\La)$, and:
\be\label{la}
{\rm(i)}\ \la(i,j)=\la(ki,kj)\quad(i,j,k\in\La),\qquad{\rm(ii)}\ |\la|:=\sum_{i\in\La}\la(0,i)=\sum_{i\in\La}\la(i,0)<\infty.
\ee

We will be interested in the continuous-time Markov process $(X_t)_{t\geq 0}$ with state space $\Si$ and generator
\bc\label{Gdef}
\dis Gf(x)&:=&\dis\sum_{i,j\in\La}\la(i,j)\big\{f\big(m^\ast_{i,j}(x)\big)-f\big(x\big)\big\}+\ga\sum_{i\in\La}\big\{f\big(m^{1\to 2}_i(x)\big)-f\big(x\big)\big\}\\[5pt]
&&\dis+\de_1\sum_{i\in\La}\big\{f\big(m^{1\to 0}_i(x)\big)-f\big(x\big)\big\}+\de'_1\sum_{i\in\La}\big\{f\big(m^{2\to 1}_i(x)\big)-f\big(x\big)\big\}\\[5pt]
&&\dis+\de_2\sum_{i\in\La}\big\{f\big(m^{\to 0}_i(x)\big)-f\big(x\big)\big\}.
\ec
Standard results \cite[Thms 4.14 and 4.29]{Swa26IPS} show that the closure of $G$ generates a Feller process $(X_t)_{t\geq 0}$, see Subsection~\ref{S:add} below for details. We call this the \emph{$(\la,\ga,\de_1,\de'_1,\de_2)$-process}. The process with $\de'_1=0$ has first been invented by Krone \cite{Kro99}, who called it the \emph{two-stage contact process}. He also considered the process with $\de_1=0$ (but possibly $\de'_1>0$), which he called the \emph{on-off contact process}.

The interpretation of the two-stage contact process in terms of young and adult individuals has already been discussed in the introduction. If the rate $\de'_1$ is positive, then adult trees grow with this rate into young trees. This seems un-biological but as Krone argued, it can be interpreted as an organism entering a state during which it cannot multiply, whence the name \emph{on-off contact process}.

Note that (\ref{la})~(i) says that the infection rates are invariant under the left action of the group on itself. For example, the group $\La$ could be the additive group~$(\Z^d,+)$. In this case~(\ref{la})~(i) says that the infection rates are translation invariant. Our setting more generally covers nearest-neighbour processes on Cayley graphs, which include regular trees, as well as long-range processes. The reason why we assume $\La$ is a group is that we need our processes to have some sort of translation-invariant structure. It seems plausible that our results hold more generally for grids that are vertex transitive and unimodular. Transitivity alone will in general not be enough, as unimodularity is essential for Lemmas \ref{L:exprate} and \ref{L:duexp} below.

\subsection{Finite and infinite survival}

For $\sig=1,2$ we define $e^\sig_i\in\Si_{\rm fin}$ by $e^\sig_i(j):=\sig$ if $i=j$ and $:=0$ otherwise. We let
\be\label{tetdef}
\tet^\sig=\tet^\sig(\la,\ga,\de_1,\de'_1,\de_2):=\P^{\txt e^\sig_0}\big[X_t\neq\un 0\ \forall t\geq 0\big]\qquad(\sig=1,2)
\ee
denote the survival probability of the $(\la,\ga,\de_1,\de'_1,\de_2)$-process started with a single young or adult tree, respectively. One has $\tet^1>0$ if and only if $\tet^2>0$, as proved in Lemma~\ref{L:tet12} below. If these survival probabilities are positive, then we say that the process exhibits \emph{finite survival} or more simply that it \emph{survives}. A simple coupling argument shows that $\tet^\sig$ is nonincreasing in $\de_1,\de'_1$, and $\de_2$, and nondecreasing in~$\ga$ (see Lemma \ref{lemmamonotonicity}).

We let $\un\sig$ denote the configuration that is constantly equal to $\sig$ $(\sig=0,1,2)$. Since the maps in (\ref{mapi}) and (\ref{mapij}) are all monotone, standard theory (compare \cite[Thm~5.13]{Swa26IPS}) tells us that there exists an invariant law $\ov\nu$ of the $(\la,\ga,\de_1,\de'_1,\de_2)$-process such that
\be
\P^{\un 2}\big[X_t\in\,\cdot\,\big]\Asto{t}\ov\nu
\ee
where $\Rightarrow$ denotes weak convergence of probability laws on $\Si=\{0,1,2\}^\La$, equipped with the product topology. Since the upper invariant law is extremal (compare \cite[Lemma~5.17]{Swa26IPS}), for a given choice of parameters either $\ov\nu=\de_{\un 0}$, or $\ov\nu$ and $\de_{\un 0}$ are singular. In the latter case we say that $\ov\nu$ is \emph{nontrivial} or that the process exhibits \emph{infinite survival}. We set
\be\label{vtetdef}
\vtet^\sig=\vtet^\sig(\la,\ga,\de_1,\de'_1,\de_2):=\int\ov\nu(\di x)1_{\txt\{x(0)>\sig\big\}}\qquad(\sig=0,1).
\ee
Then under the upper invariant law, $\vtet^0$ is the density of trees (young or old) while $\vtet^1$ is the density of adult trees. In Lemma \ref{L:tet12} below we prove that $\vtet^0>0$ if and only if $\vtet^1>0$.

For the contact process, it is well-known that finite survival is equivalent to infinite survival. This has a one-line proof based on the well-known self-duality of the contact process. Krone discovered a duality between two-stage contact processes and on-off contact processes that allows one to prove, in one line, that
\be\label{fininf}
\vtet^\sig(\la,\ga,\de_1,\de'_1,\de_2)=\theta^{2-\sig}(\la^\dgg,\ga,\de'_1,\de_1,\de_2),
\ee
where $\la^\dgg(i,j):=\la(j,i)$ denote the \emph{reversed infection rates}. Note that on the right-hand side of~(\ref{fininf}), not only is $\la$ replaced by $\la^\dgg$, but also the roles of $\de_1$ and $\de'_1$ are reversed. See formula~(\ref{tetvar}) below for details. In particular,~(\ref{fininf}) implies that finite survival of the ~$(\la,\ga,\de_1,\de'_1,\de_2)$-process is equivalent to infinite survival of the $(\la^\dgg,\ga,\de'_1,\de_1,\de_2)$-process. Since the two-stage contact process is not self-dual (except when $\de_1=0$), this does not answer Krone's question whether for a given process, finite and infinite survival are equivalent.

To study this question, we let $D:=\big\{(\ga,\de_1,\de'_1,\de_2):\ga>0,\ \de_1,\de'_1,\de_2\geq 0\big\}$ denote our parameter space. We fix $\la$ and define
\bc\label{Surv}
\dis D_{\rm fin}&:=&\dis\big\{(\ga,\de_1,\de'_1,\de_2)\in D:\tet^1(\la,\ga,\de_1,\de'_1,\de_2)>0\big\},\\[5pt]
\dis D_{\rm inf}&:=&\dis\big\{(\ga,\de_1,\de'_1,\de_2)\in D:\vtet^0(\la,\ga,\de_1,\de'_1,\de_2)>0\big\}.
\ec
For finite range processes on $\Z^d$, the results of \cite{BG90,BG94} imply that $D_{\rm fin}$ is an open subset of $D$, that is, critical processes die out. Although this is not explicitly mentioned in these papers, Foxall \cite[Section~3.5]{Fox15} observed that the proof methods of \cite{BG90,BG94} also allow one to conclude that $D_{\rm fin}\sub D_{\rm inf}$. By applying this both to the two-stage contact process and to its dual, the on-off contact process, using (\ref{fininf}), he obtained that $D_{\rm fin}=D_{\rm inf}$, solving Krone's question for finite range processes on $\Z^d$.

The methods of \cite{BG90,BG94} are restricted to $\Z^d$. For contact processes on general groups, it is not known if the critical process dies out. We will show, however, that the critical points for finite and infinite survival of the two-stage contact process on general groups are equal. More precisely, we will prove the following theorem, which is the first main result of the paper. Below, we let $\ring A$ denote the interior of a set $A\sub D$.

\bt[Finite and infinite survival]
The\label{T:fininf} sets defined in (\ref{Surv}) satisfy $\ring D_{\rm fin}=\ring D_{\rm inf}$.
\et

\subsection{Sharpness of the phase transition}\label{S:sharp}

For any $x\in\Si_{\rm fin}$, we let $|x|:=\sum_{i\in\La}1_{\{x(i)>0\}}$ denote the number of trees (young or adult). We call the constant $r(\la,\ga,\de_1,\de'_1,\de_2)$ from the following lemma the \emph{exponential growth rate} of the $(\la,\ga,\de_1,\de'_1,\de_2)$-process. 
\bl[The exponential growth rate]
The\label{L:exprate} $(\la,\ga,\de_1,\de'_1,\de_2)$-process satisfies
\be\label{exprate}
\lim_{t\to\infty}\frac{1}{t}\log\E^{e^2_0}\big[|X_t|]=\inf_{t>0}\frac{1}{t}\log\E^{e^2_0}\big[|X_t|]=:r
\ee
for some constant $r=r(\la,\ga,\de_1,\de'_1,\de_2)\in\R$. The quantity $\E^{e^2_0}\big[|X_t|]$ is the same for the $(\la,\ga,\de_1,\de'_1,\de_2)$-process and the $(\la^\dgg,\ga,\de'_1,\de_1,\de_2)$-process, so as a consequence
\be\label{durat}
r(\la,\ga,\de_1,\de'_1,\de_2)=r(\la^\dgg,\ga,\de'_1,\de_1,\de_2).
\ee
\el

We fix $\la$ and set
\be\label{Dpos}
D_{\rm pos}:=\big\{(\ga,\de_1,\de'_1,\de_2)\in D:r(\la,\ga,\de_1,\de'_1,\de_2)\geq 0\big\}.
\ee
We will derive Theorem~\ref{T:fininf} as a simple consequence of formula (\ref{fininf}), Lemma~\ref{L:exprate}, and the following theorem, that is our second main result. In words, the theorem says that subcritical $(\la,\ga,\de_1,\de'_1,\de_2)$-processes die out exponentially fast. 

\bt[Sharpness of the phase transition]
The\label{T:sharp} sets defined in (\ref{Surv}) and (\ref{Dpos}) satisfy $\ring D_{\rm fin}=\ring D_{\rm pos}$.
\et

The proof of Lemma~\ref{L:exprate} is deferred to Section~\ref{S:exp} and the proof of Theorem~\ref{T:sharp} is deferred to Section~\ref{sec:proofTsharp}. We now prove Theorem~\ref{T:fininf}.\med

\bpro[of Theorem~\ref{T:fininf}]
Let us write $D_{\rm fin}(\la)$, $D_{\rm inf}(\la)$, and $D_{\rm pos}(\la)$ to stress the dependence of these sets on the infection rates $\la$, and for any set $A\sub D$, let us write
\be
\widehat A:=\big\{(\ga,\de'_1,\de_1,\de_2):(\ga,\de_1,\de'_1,\de_2)\in A\big\},
\ee
where we have flipped the second and third coordinates. Formula (\ref{fininf}) tells us that $D_{\rm inf}(\la)={\widehat D}_{\rm fin}(\la^\dgg)$ while formula (\ref{durat}) of Lemma~\ref{L:exprate} tells us that $D_{\rm pos}(\la)=\widehat D_{\rm pos}(\la^\dgg)$. Using also Theorem~\ref{T:sharp}, this allows us to conclude that
\be
\ring D_{\rm fin}(\la)=\ring D_{\rm pos}(\la)=\ring{\widehat D}_{\rm pos}(\la^\dgg)=\ring{\widehat D}_{\rm fin}(\la^\dgg)=\ring D_{\rm inf}(\la).
\ee
\epro

\subsection{Discussion and outline}

We have just seen that Theorem~\ref{T:fininf} is a simple consequence of formula (\ref{fininf}), Lemma~\ref{L:exprate}, and Theorem~\ref{T:sharp}. In formula (\ref{tetvar}), Lemma~\ref{L:duexp}, and Theorem~\ref{T:sharp2} below, we will prove more general versions of these results for additive interacting particle systems whose local state space can be any finite lattice. We note that if this lattice is distributive, then such a particle system has a percolation representation in terms of open paths, as proved in \cite{Swa26perc}. For our general results, we do not need distributivity, however.

Formula (\ref{tetvar}) and Lemma~\ref{L:duexp} follow from simple duality arguments but the proof of Theorem~\ref{T:sharp2}, which establishes sharpness of the phase transition for a general class of additive interacting particle systems, requires more work. We use the technique of \cite{Swa18}, where this was proved for contact processes. Very briefly, the idea of the proof is as follows. By definition, an \emph{eigenmeasure} is a locally finite measure $\mu$ on the space of nonzero configurations such that if the process is started in the initial ``law'' $\mu$ (which may be an infinite measure), then its ``law'' at time $t$ is $e^{\la t}\mu$. Here $\la$ is called the \emph{eigenvalue}, see formula (\ref{eig}) below. Generalising results from \cite{Swa09}, in Proposition~\ref{P:eigmeas} below we show that additive interacting particle systems have an eigenmeasure $\mu$ whose eigenvalue is the exponential growth rate. Via duality, such an eigenmeasure gives rise to an eigenfunction $h$ of the generator of the dual process, see Lemma~\ref{L:lemmah} below. In particular, if the exponential growth rate is zero, then $h$ is a harmonic function. Theorem~\ref{T:sharp2} is then proved by showing that by changing the rates a little bit, it is possible to transform $h$ into a subharmonic function that can be used to prove that finite processes survive.

We note that for contact processes, it is established in \cite[Thm~1.2~(b)]{Swa09} that the exponential growth rate is Lipschitz continuous as a function of the death rate. We had originally planned to establish an analogous result for a general class of additive interacting particle systems, and in particular, we had hoped to prove that the exponential growth rate of the two-stage contact process depends continuously on its parameters. This turned out to be harder than expected and is still open.

\section{Preliminaries}

\subsection{Additive maps and their duals}

We will need some general theory for additive interacting particle systems whose local state space can be any finite lattice. We recall that a \emph{lattice} is a partially ordered set $\Lt$ in which each two elements $a,b\in\Lt$ have a greatest lower bound $a\wedge b$ and least upper bound $a\vee b$, called the \emph{infimum} and \emph{supremum} of $a$ and $b$. Each finite lattice $\Lt$ contains a least element $0$ and greatest element $\top$. We fix a countable set $\La$ and let $\Si=\Si(\La):=\Lt^\La$ denote the space of all functions $x\cn\La\to\Lt$. We let $\Si_{\rm fin}=\Si_{\rm fin}(\La)$ denote the subspace consisting of all $x\in\Si$ for which $\{i\in\La:x(i)\neq 0\}$ is a finite set. We equip $\Si=\Lt^\La$ with the product topology and the product order.

The space $\Si$ is a \emph{complete} lattice, meaning that each subset $A\sub\Si$ has a greatest lower bound $\inf A$ and least upper bound $\sup A$. For a sequence $(x_n)_{n\in\N}$ in $\Si$ we adopt the notation~$\bigwedge_nx_n:=\inf\{x_n:n\in\N\}$ and $\bigvee_nx_n:=\sup\{x_n:n\in\N\}$. These are the pointwise infimum and supremum of the functions $x_n$. The space $\Si$ contains a least element $\un 0$ and greatest element $\un\top$, which are the functions that are constantly equal to $0$ and $\top$, respectively.

A \emph{dual} of a partially ordered set $\Lt$ is a partially ordered set $\Lt'=(\Lt',\leq)$ together with a bijection $\Lt\ni a\mapsto a'\in\Lt'$ such that
\be
a\leq b\quad\mbox{if and only if}\quad a'\geq b'.
\ee
All duals of $\Lt$ are naturally isomorphic to each other, and setting $a'':=a$, we can naturally identify $\Lt$ with the dual of $\Lt'$. For example, the totally ordered set $\Lt=\{0,1,2\}$ is dual to itself via the map $a\mapsto a':=2-a$ $(a\in\Lt)$. We fix a dual $\Lt'$ of $\Lt$ and define $\Si'$ and $\Si'_{\rm fin}$ as before with $\Lt$ replaced by $\Lt'$. Then $\Si'$ is dual to $\Si$ via the map $x\mapsto x'$ that is defined pointwise as $x'(i):=x(i)'$ $(i\in\La)$.

A map $m\cn\Si\to\Si$ is \emph{additive} if 
\be
m(\un 0)=\un 0\quand m(x\vee y)=m(x)\vee m(y)\qquad(x,y\in\Si),
\ee
and \emph{countably additive} if moreover $m\big(\bigvee_nx_n\big)=\bigvee_nm(x_n)$ for each sequence $x_n$ in $\Si$. It is then easy to see that each continuous additive map is also countably additive. 
We define $\phi\cn\Si\times\Si'\to\{0,1\}$ by
\be\label{phidef}
\phi(x,y):=1_{\txt\{x\leq y'\}}=1_{\txt\{x'\geq y\}}
\qquad(x\in\Si,\ y\in\Si').
\ee
If $m\cn\Si\to\Si$ and $\hat m\cn\Si'\to\Si'$ are maps, then we say that $\hat m$ is \emph{dual} to $m$ with respect to the duality function $\phi$ if
\be\label{mdual}
\phi\big(m(x),y\big)=\phi\big(x,\hat m(y)\big)
\qquad(x\in\Si,\ y\in\Si').
\ee
We say that a map $m\cn\Si\to\Si$ is \emph{local} if there exists a finite set $\De\sub\La$ and a map $n\cn\Si(\De)\to\Si(\De)$ such that $m(x)(i)=n(x\!\restriction_\De)(i)$ $(i\in\De)$ and $m(x)(i)=x(i)$ $(i\in\La\beh\De)$, where $x\!\restriction_\De$ denotes the restriction of $x$ to $\De$. We cite the following lemma from \cite[Lemma~1]{Swa26perc}.

\bl[Dual maps]
Let\label{L:dumaps} $m\cn\Si\to\Si$ be countably additive. Then $m$ has a unique dual map $\hat m\cn\Si'\to\Si'$ with respect to the duality function $\phi$ from (\ref{phidef}). The dual map $\hat m$ is also countably additive. If $m$ is a local map, then so is $\hat m$.
\el

It is sometimes more convenient to work with
\be
\psi(x,y):=1-\phi(x,y)=1_{\{x\not\leq y'\}}\qquad(x\in\Si,\ y\in\Si').
\ee
Clearly, if two maps $m$ and $\hat m$ are dual with respect to $\phi$ in the sense defined in (\ref{mdual}), then they are also dual with respect to $\psi$. 

We already observed that $\Lt=\{0,1,2\}$ is dual to itself via the map $a\mapsto a':=2-a$ $(a\in\Lt)$. In this case, the duality function in (\ref{phidef}) takes the form
\be \label{dualityfunction}
\phi(x,y):=1_{\txt\{x(i)+y(i)\leq 2\ \forall i\in\La\}}\qquad(x,y\in\Si).
\ee
It is straightforward to check (and has already been done by Krone in \cite{Kro99}) that the local maps in (\ref{mapi}) and (\ref{mapij}) are additive and their duals are given by
\be\ba{r@{\,}c@{\,}l@{\quad}r@{\,}c@{\,}l}\label{dumaps}
\dis\widehat{m^{1\to 2}_i}&=&\dis m^{1\to 2}_i,&
\dis\widehat{m^{1\to 0}_i}&=&\dis m^{2\to 1}_i\\[5pt]
\dis\widehat{m^{2\to 1}_i}&=&\dis m^{1\to 0}_i,&
\dis\widehat{m^{\to 0}_i}&=&\dis m^{\to 0}_i\\[5pt]
\multicolumn{6}{c}{\widehat{m^\ast_{i,j}}=m^\ast_{j,i}.}
\ec
A nice graphical explanation of these duality relations can be found in \cite[Fig~1]{Swa26perc}.

For $a\in\Lt$ and $i\in\La$ we define $e^a_i\in\Si_{\rm fin}$ by $e^a_i(j):=a$ if $i=j$ and $:=0$ otherwise. If~$m\cn\Si\to\Si$ is countably additive, then for each $i,j\in\La$ we define $m[i,j]\cn\Lt\to\Lt$ by
\be\label{mat}
m[i,j](a):=m(e^a_i)(j)\qquad(a\in\Lt,\ i,j\in\La).
\ee
Since $m$ is countably additive,
\be\label{mij}
m(x)(j)=\bigvee_{i\in\La}m[i,j]\big(x(i)\big).
\ee
We use similar notation for the dual map $\hat m$. It is then straightforward to check that
\be\label{mijdu}
\hat m[i,j]=\widehat{m[j,i]}\qquad(i,j\in\La),
\ee
where $\widehat{m[j,i]}$ is the dual of the map $m[i,j]\cn\Lt\to\Lt$.

\subsection{Additive particle systems}\label{S:add}

Generalising the set-up of (\ref{Gdef}), we will be interested in interacting particle systems $(X_t)_{t\geq 0}$ with state space $\Si=\Si(\La)$ and generator of the form
\be\label{Gf}
Gf(x):=\sum_{m\in\Gi}r_m\big\{f\big(m(x)\big)-f\big(x\big)\big\},
\ee
where $\Gi$ is a set whose elements are additive local maps $m\cn\Si\to\Si$ and $(r_m)_{m\in\Gi}$ are nonnegative rates. In this context, $\La$ is called the \emph{grid}\footnote{This is often called the \emph{lattice}, but we will need that word in its order-theoretic meaning.} and $\Lt$ is called the \emph{local state space}. We call~(\ref{Gf}) a \emph{random mapping representation} of the generator $G$. We let ${\rm id}\cn\Lt\to\Lt$ denote the identity map and we let $o\cn\Lt\to\Lt$ denote the zero map defined as $o(a):=0$ $(a\in\Lt)$. We use similar notation for the identity map and zero map on $\Lt'$. It is easy to see that $\hat{\rm id}={\rm id}$ and $\hat o=o$, where the maps on the left-hand side are defined on $\Lt'$. We make the following assumptions on the rates:
\be\ba{l}\label{sumad}
\dis{\rm(i)}\ \sup_{i\in\La}\sum_{m\in\Gi}r_m1_{\txt\{m[i,i]\neq{\rm id}\}}<\infty,\\[5pt]
\dis{\rm(ii)}\ \sup_{i\in\La}\sum_{m\in\Gi}r_m\sum_{j\in\La\beh\{i\}}1_{\txt\{m[j,i]\neq o\}}<\infty,\\[5pt]
\dis{\rm(iii)}\ \sup_{i\in\La}\sum_{m\in\Gi}r_m\sum_{j\in\La\beh\{i\}}1_{\txt\{m[i,j]\neq o\}}<\infty.
\ec
By \cite[Lemma~9]{Swa26perc} and \cite[Thms 4.14 and 4.29]{Swa26IPS}, these conditions guarantee that the closure of $G$, which is initially defined only for functions depending on finitely many coordinates, generates a Markov process with state space $\Si$. This Markov process can moreover be constructed from a graphical representation, as we now explain.

Let $\om$ be a Poisson point set on $\Gi\times\R$ with intensity $\rho$ given by
\be\label{rho}
\rho\big(\{m\}\times[s,t]\big):=r_m(t-s)\qquad(m\in\Gi,\ s\leq t).
\ee
Recall that a function $f$, defined on a real interval $I$ and taking values in a metrisable space, is called \emph{cadlag} (from the French continue \`a droite, limite \`a gauche) if it is right-continuous with left limits $f(t-):=\lim_{s\up t}f(s)$ in each point $t\in I$ that can be approximated from the left. By \cite[Lemma~9]{Swa26perc} and \cite[Thm~4.14]{Swa26IPS}, almost surely, for each $x\in\Si$ and $s\in\R$, there exists a unique cadlag function $X^{s,x}\cn[s,\infty)\to\Si$ with initial state $X^{s,x}_s=x$ that satisfies the evolution equation
\be\label{evol}
X^{s,x}_t:=\left\{\ba{ll}
\dis m(X^{s,x}_{t-})\quad&\mbox{if }(m,t)\in\om,\\[5pt]
X^{s,x}_{t-}\quad&\mbox{otherwise,}\ea\right.\qquad(s<t).
\ee
Moreover, setting $\Xb_{s,t}(x):=X^{s,x}_t$ $(s\leq t,\ x\in\Si)$ defines a collection $(\Xb_{s,t})_{s\leq t}$ of continuous maps from $\Si$ into itself such that
\be\label{oldflow}
\Xb_{s,s}=1\quand\Xb_{t,u}\circ\Xb_{s,t}=\Xb_{s,u}\qquad(s\leq t\leq u),
\ee
where $1$ denotes the identity map. We call $(\Xb_{s,t})_{s\leq t}$ the \emph{stochastic flow} defined by the \emph{graphical representation} $\om$. Since restrictions of a Poisson point process to disjoint sets are independent, the stochastic flow $(\Xb_{s,t})_{s\leq t}$ has \emph{independent increments} in the sense that for each $t_0<\cdots<t_n$, the random variables $(\Xb_{t_{k-1},t_k})_{1\leq k\leq n}$ are independent. It is also \emph{stationary} in the sense that $(\Xb_{s,t})_{s\leq t}$ is equally distributed with $(\Xb_{r+s,r+t})_{s\leq t}$ for each $r\in\R$.

Similarly, it can be shown that for each $u\in\R$ and $y\in\Si'$, there exists a unique cadlag function $Y^{u,y}\cn(-\infty,u]\to\Si'$ with final state $Y^{u,y}_u=y$ that satisfies the evolution equation
\be\label{duevol}
Y^{u,y}_{t-}:=\left\{\ba{ll}
\dis\hat m(Y^{u,y}_t)\quad&\mbox{if }(m,t)\in\om,\\[5pt]
Y^{u,y}_t\quad&\mbox{otherwise,}\ea\right.\qquad(t\leq u).
\ee
Setting $\Yb_{u,t}(y):=Y^{u,y}_t$ $(u\geq t,\ y\in\Si')$ then defines a collection $(\Yb_{u,t})_{u\geq t}$ of continuous maps from $\Si'$ into itself such that
\be\label{duflow}
\Yb_{s,s}=1\quand\Yb_{t,s}\circ\Yb_{u,t}=\Yb_{u,s}\qquad(u\geq t\geq s).
\ee
We call $(\Yb_{u,t})_{u\geq t}$ the \emph{backward stochastic flow} defined by the \emph{dual graphical representation} $\hat\om:=\{(\hat m,t):(m,t)\in\om\}$. If $u\in\R$ and $Y_0$ is an $\Si'$-valued random variable with law $\mu$, independent of $\om$, then setting
\be\label{YY}
Y_t:=\Yb_{u,u-t}(Y_0)\qquad(t\geq 0)
\ee
defines a Feller process $(Y_t)_{t\geq 0}$ with generator
\be\label{hatGf}
\hat Gf(y):=\sum_{m\in\Gi}r_m\big\{f\big(\hat m(y)\big)-f\big(y\big)\big\}.
\ee
We note that by (\ref{mijdu}) and the facts that $\hat{\rm id}={\rm id}$ and $\hat o=o$, the rates of this dual process also satisfy the condition (\ref{sumad}). Because of time reversal, somewhat unusually, the sample paths of the Markov process defined in (\ref{YY}) are \emph{caglad}, that is, left-continuous with right limits. By (\ref{dumaps}), if $G$ is the generator of the $(\la,\ga,\de_1,\de'_1,\de_2)$-process, written in the form (\ref{Gdef}), then $\hat G$ is the generator of the $(\la^\dgg,\ga,\de'_1,\de_1,\de_2)$-process. By \cite[Thm~2]{Swa26perc}, almost surely $\Xb_{s,t}$ and $\Yb_{t,s}$ are countably additive for each $s\leq t$ and
\be\label{flowdu}
\phi\big(\Xb_{s,t}(x),y\big)=\phi\big(x,\Yb_{t,s}(y)\big)
\qquad(x\in\Si,\ y\in\Si'),
\ee
that is, $\Yb_{t,s}=\hat\Xb_{s,t}$ is the dual of $\Xb_{s,t}$. In particular, this implies that the Markov processes $(X_t)_{t\geq 0}$ and $(Y_t)_{t\geq 0}$ are dual in the sense that
\be\label{dual}
\E^x\big[\phi(X_t,y)\big]=\E^y\big[\phi(x,Y_t)\big]\qquad(t\geq 0,\ x\in\Si,\ y\in\Si'),
\ee
where $\P^x$ denotes the law of the process $(X_t)_{t\geq 0}$ started in $X_0=x$, $\E^x$ denotes expectation with respect to $\P^x$, and similar notation applies to the dual process $(Y_t)_{t\geq 0}$. By (\ref{dualityfunction}), it follows that the $(\la,\ga,\de_1,\de'_1,\de_2)$-process $(X_t)_{t\geq 0}$ and the $(\la^\dgg,\ga,\de'_1,\de_1,\de_2)$-process $(Y_t)_{t\geq 0}$ satisfy
\be\label{ourdual}
\P^x\big[X_t(i)+y(i)\leq 2\ \forall i\in\La\big]=
\P^y\big[x(i)+Y_t(i)\leq 2\ \forall i\in\La\big]
\qquad(t\geq 0,\ x,y\in\Si).
\ee

\subsection{Functions of polynomial growth}

By \cite[Lemma~9]{Swa26perc} and \cite[Prop~4.24]{Swa26IPS}, $\Xb_{s,t}$ and $\Yb_{t,s}$ map the space $\Si_{\rm fin}$ into itself, almost surely for all $s\leq t$. Since the set $\Si_{\rm fin}$ is countable, the Markov process $(X_t)_{t\geq 0}$ started in an initial condition $X_0=x\in\Si_{\rm fin}$ is just a continuous-time Markov chain, and similarly for the dual process $(Y_t)_{t\geq 0}$. We let $P_t(x_1,x_2)$ denote the transition probabilities of $(X_t)_{t\geq 0}$ and denote its generator (also called Q-matrix) by $G$. This generator is given by
\be\label{matrix}
G(x_1,x_2):=\sum_{m\in\Gi}r_m1_{\txt\{m(x_1)=x_2\}}\quad(x_1\neq x_2),
\ee
with the diagonal entries determined by the condition that $\sum_{x_2\in\Si_{\rm fin}}G(x_1,x_2)=0$ $(x_1\in\Si_{\rm fin})$. Similarly, we let $\hat P_t(y_1,y_2)$ denote the transition probabilities of $(Y_t)_{t\geq 0}$ and we denote its generator by $\hat G$. For $x\in\Si$ we set $|x|:=\sum_{i\in\La}1_{\{x(i)\neq 0\}}$, so that $\Si_{\rm fin}=\{x\in\Si:|x|<\infty\}$. Similar notation applies to $\Si'$. By \cite[Prop~4.24]{Swa26IPS} and \cite[Lemma~8]{Swa26perc}, one has the bounds
\be\label{expbd}
\left.\bac
\dis\E^x\big[|X_t|\big]&\leq&\dis|x|\ex{K_\up t},\\[5pt]
\dis\E^y\big[|Y_t|\big]&\leq&\dis|y|\ex{K_\down t},
\ea\right\}
\qquad(t\geq 0,\ x\in\Si_{\rm fin},\ y\in\Si'_{\rm fin}),
\ee
where
\be\label{Kup}
K_\down:=\sup_{i\in\La}\sum_{m\in\Gi}r_m\big(\sum_{j\in\La}1_{\{m[j,i]\neq o\}}-1\big)
\quand
K_\up:=\sup_{i\in\La}\sum_{m\in\Gi}r_m\big(\sum_{j\in\La}1_{\{m[i,j]\neq o\}}-1\big)
\ee
are finite by (\ref{sumad}). Formula (\ref{expbd}) shows that $P_tf$ is well-defined provided $f$ grows at most linearly in $|x|$.

We will need some technical results that allow us to do generator calculations for functions that grow faster than linearly. To this aim, we introduce the additional assumption that there exists an $1\leq R<\infty$ such that
\be\label{Rbd}
\big|\big\{(i,j)\in\La^2:m[i,j]\neq{\rm Id}[i,j]\big\}\big|\leq R\qquad(m\in\Gi),
\ee
where ${\rm Id}\cn\Si\to\Si$ denotes the identity map and consequently ${\rm Id}[i,i]={\rm id}$ and ${\rm Id}[i,j]=o$ $(i\neq j)$. In view of (\ref{mijdu}) and the facts that $\hat{\rm id}={\rm id}$ and $\hat o=o$, if (\ref{Rbd}) is satisfied for the process $(X_t)_{t\geq 0}$, then also for the dual process $(Y_t)_{t\geq 0}$. For $(\ga,\de_1,\de'_1,\de_2)$-processes, (\ref{Rbd}) is satisfied with $R=1$. We define a metric $d$ on $\Si_{\rm fin}$ by
\be\label{ddef}
d(x,y):=\sum_{i\in\La}1_{\{x(i)\neq y(i)\}}\qquad(x,y\in\Si_{\rm fin}).
\ee
Then $d(x,\un 0)=|x|$ $(x\in\Si_{\rm fin})$. We observe that $m(x)(j)\neq x(j)$ implies that there exists an $i\in\La$ with $m[i,j]\neq{\rm Id}(i,j)$, so (\ref{Rbd}) implies that
\be\label{Rbd2}
d\big(m(x),x\big)\leq R\qquad(m\in\Gi).
\ee

For $k\in\N$, we let $\Fi_k=\Fi_k(\Si_{\rm fin})$ denote the space of all functions $f\cn\Si_{\rm fin}\to\R$ that are of \emph{$k$-th order polynomial growth} in the sense that there exist a constant $C<\infty$ such that
\be\label{Fkdef}
|f(x)|\leq C\big(|x|^k+1\big)\qquad(x\in\Si_{\rm fin}).
\ee
In particular, $\Fi_0$ are the bounded functions, while functions in $\Fi_1$ and $\Fi_2$ are of \emph{linear growth} and \emph{quadratic growth}, respectively. We let $\Fi_\infty=\Fi_\infty(\Si_{\rm fin}):=\bigcup_{k=0}^\infty\Fi_k(\Si_{\rm fin})$ denote all functions of \emph{polynomial growth} (or arbitrary order).

We will say that a sequence of functions $f_n\cn S\to\R$ converges to a limit $f$ in a \emph{polynomially bounded pointwise} way, or \emph{in pbp-sense} for short, if $f_n(x)\to f(x)$ for all $x\in S$ and there exist constants $k\in\N$ and $C<\infty$ such that $|f_n(x)|\leq C(|x|^k+1)$ for all $n$. We also write $\pbplim_{n\to\infty}f_n=f$. Polynomially bounded pointwise convergence for functions $f_t$ depending on a real parameter $t$ is defined similarly. The following theorem is the main result of this section.

\bt[Functions of polynomial growth]
Assume\label{T:poly} (\ref{Rbd}). Then the operators $G$ and $P_t$ $(t\geq 0)$ map the space $\Fi_\infty$ into itself. Moreover $\pbplim_{n\to\infty}f_n=f$ implies $\pbplim_{n\to\infty}P_tf_n=P_tf$ for each $t\geq 0$. Finally,
\be\label{pbplim}
\pbplim_{t\down 0}\,t^{-1}(P_t-1)f=Gf\qquad(f\in\Fi_\infty).
\ee
\et

Using Theorem~\ref{T:poly}, it is easy to prove (see Lemma~\ref{L:rdif} in the appendix) that for each $f\in\Fi_\infty$, the functions $t\mapsto P_tf(x)$ are continuously differentiable for all $x\in S$, and
\be\label{poly}
\dif{t}P_tf=GP_tf=P_tGf\qquad(t\geq 0).
\ee
Before we give the proof of Theorem~\ref{T:poly}, to demonstrate its usefulness, we state one more consequence of it that will be needed in what follows.

We view a real function $f$ on $\Si_{\rm fin}$ as a column vector and a probability law $\mu$ as a row vector and consequently write $\mu P_t(x_2):=\sum_{x_1\in\Si_{\rm fin}}\mu(x_1)P_t(x_1,x_2)$ for the law at time $t$ of the process started in the initial law $\mu$, and $\mu f:=\sum_{x\in\Si_{\rm fin}}\mu(x)f(x)$ for the expectation of $f$ under $\mu$. We let $\cov_\mu(f,g):=\mu(fg)-(\mu f)(\mu g)$ denote the covariance of $f$ and $g$ under $\mu$. We note that $fg\in\Fi_\infty$ for all $f,g\in\Fi_\infty$, which allows us to define
\be
\Ga(f,g):=G(fg)-(Gf)g-f(Gg)\qquad(f,g\in\Fi_\infty).
\ee
One can check that
\bc \label{eqgamma}
\dis\Ga(f,g)(x_1)
&=&\dis\sum_{x_2 \in \Si_{\rm fin}}G(x_1,x_2)\big(f(x_2)-f(x_1)\big)\big(g(x_2)-g(x_1)\big),\\[5pt]
&=&\dis\sum_{m\in\Gi}r_m\big(f\big(m(x_1)\big)-f\big(x_1\big)\big)\big(g\big(m(x_1)\big)-g\big(x_1\big)\big).
\ec

\bl[Covariance formula]
Assume\label{L:cov} (\ref{Rbd}) and let $\mu$ be a probability measure on $\Si_{\rm fin}$ such that $\sum_{x\in\Si_{\rm fin}}\mu(x)|x|^k<\infty$ for all $k\in\N$. Then the same holds for the measures $\mu P_t$ $(t\geq 0)$, and for each $f,g\in\Fi_\infty$, one has
\be
\cov_{\mu P_t}(f,g)=\cov_\mu(P_tf,P_tg)+\int_0^t\mu P_{t-s}\Ga(P_sf,P_sg)\,\di s\qquad(t\geq 0).
\ee
\el

\bpro
The fact that $\sum_{x\in\Si_{\rm fin}}\mu P_t(x)|x|^k<\infty$ for all $k\in\N$ and $t\geq 0$ follows from the fact, proved in Theorem~\ref{T:poly}, that $P_t$ maps the space $\Fi_\infty$ into itself. The covariance formula now follows from a standard calculation, which can be found in \cite[Prop~2.3]{Swa26IPS}, provided we show that
\bc
\dis\dif{t_1}P_{t_1}\big((P_{t_2}f)(P_{t_3}g)\big)&=&\dis P_{t_1}G\big((P_{t_2}f)(P_{t_3}g)\big),\\[5pt]
\dis\dif{t_2}P_{t_1}\big((P_{t_2}f)(P_{t_3}g)\big)&=&\dis P_{t_1}\big((GP_{t_2}f)(P_{t_3}g)\big),\\[5pt]
\dis\dif{t_3}P_{t_1}\big((P_{t_2}f)(P_{t_3}g)\big)&=&\dis P_{t_1}\big((P_{t_2}f)(GP_{t_3}g)\big).
\ec
In our present setting, these formulas follow from Theorem~\ref{T:poly}, the properties of polynomially bounded pointwise convergence, and Lemma~\ref{L:rdif} in the appendix.
\epro

We next start to prepare for the proof of Theorem~\ref{T:poly}.

\bl[The rate of jumping from a state]
There\label{L:Gdia} exists a constant $K<\infty$ such that
\be
-G(x,x)\leq K|x|\qquad(x\in\Si_{\rm fin}).
\ee
\el

\bpro
Let
\bc\label{K01}
\dis K_0&:=&\dis\sup_{i\in\La}\sum_{m\in\Gi}r_m1_{\txt\{m[i,i]\neq{\rm id}\}}<\infty,\\[5pt]
\dis K_1&:=&\dis\dis\sup_{i\in\La}\sum_{m\in\Gi}r_m\sum_{j\in\La\beh\{i\}}1_{\txt\{m[i,j]\neq o\}}<\infty,
\ec
which are finite by (\ref{sumad}) (i) and (iii). In view of (\ref{mij}), for each $m\in\Gi$ and $x\in\Si_{\rm fin}$, if $m(x)\neq x$, then there must exist $i,j\in\La$ such that $x(i)\neq 0$ and either $i=j$ and $m[i,i]\neq{\rm id}$ or $i\neq j$ and $m[i,j]\neq o$. It follows that
\be\ba{l}
\dis-G(x,x)=\sum_{m\in\Gi}r_m1_{\txt\{m(x)\neq x\}}\\[5pt]
\dis\quad\leq\sum_{m\in\Gi}r_m\sum_{i\in\La}1_{\txt\{x(i)\neq 0\}}\Big(1_{\txt\{m[i,i]\neq{\rm id}\}}+\sum_{j\in\La\beh\{i\}}1_{\txt\{m[i,j]\neq o\}}\Big)\\[5pt]
\dis\quad\leq(K_0+K_1)\sum_{i\in\La}1_{\txt\{x(i)\neq 0\}}=(K_0+K_1)|x|.
\ec
\epro

\bl[Bound on the generator]
Assume\label{L:GFk} (\ref{Rbd}). Then for each $k\geq 0$, the operator $G$ maps $\Fi_k$ into $\Fi_{k+1}$.
\el

\bpro
Assume that $|f(x)|\leq C\big(|x|^k+1\big)$ $(x\in\Si_{\rm fin})$ for some $C<\infty$. Then using (\ref{Rbd2}), we can estimate
\bc
\dis|Gf(x)|
&=&\dis\big|\sum_{m\in\Gi}r_m\big\{f\big(m(x)\big)-f\big(x\big)\big\}\big|\\[5pt]
&\leq&\dis\sum_{m\in\Gi}r_m1_{\txt\{m(x)\neq x\}}\big\{C\big((|x|+R)^k+1\big)+C\big(|x|^k+1\big)\big\}\\
&\leq&\dis CK|x|\big\{(|x|+R)^k+|x|^k+2\big\},
\ec
where $K$ is the constant from Lemma~\ref{L:Gdia}.
\epro

We will derive Theorem~\ref{T:poly} from the abstract Theorem~\ref{T:Lbp} in the appendix. To this aim, we need a sequence of Lyapunov functions. For $k\geq 1$, we define $g_k\cn\N\to\N$ by
\be
g_k(n):=\prod_{i=0}^{k-1}\big(n+iR\big)\qquad(n\geq 0,\ k\geq 1),
\ee
where $R$ is the constant from (\ref{Rbd}) and we set $g_0(n):=1$ for all $n\geq 0$. Then, for any $k\geq 1$,
\be
g_k(n+R)-g_k(n)=\big((n+kR)-n\big)\prod_{i=1}^{k-1}(n+iR)=kRg_{k-1}(n+R)
\ee
and hence
\be\label{diftrick}
n\big(g_k(n+R)-g_k(n)\big)=kRng_{k-1}(n+R)=kRg_k(n).
\ee
We set
\be
L_k(x):=g_k(|x|)+1\qquad(k\geq 0,\ x\in\Si_{\rm fin}).
\ee

\bl[Lyapunov functions]
Assume\label{L:Lya} (\ref{Rbd}). Then for each $k\geq 0$, there exists a constant $\la_k\in\R$ such that $GL_k\leq\la_k L_k$.
\el

\bpro
For each $m\in\Gi$, if $L_k(m(x))>L_k(x)$, there must exist $i,j\in\La$ with $i\neq j$ such that $x(i)\neq 0$ and $m[i,j]\neq o$. Together with (\ref{Rbd2}), this allows us to estimate
\bc
\dis GL_k(x)
&=&\dis\sum_{m\in\Gi}r_m\big\{L_k\big(m(x)\big)-L_k\big(x\big)\big\}\\[5pt]
&\leq&\dis\sum_{m\in\Gi}r_m1_{\txt\{\exists i\neq j\mbox{ s.t.\ }x(i)\neq 0, m[i,j]\neq o\}}\big\{g_k(|x|+R)-g_k(|x|)\big\}\\[5pt]
&\leq&\dis\sum_{m\in\Gi}r_m\sum_{i\in \La}1_{\txt\{x(i)\neq 0\}}\sum_{j\in \La\backslash{\{i\}}}1_{\txt\{m[i,j]\neq o\}}\big\{g_k(|x|+R)-g_k(|x|)\big\}\\[5pt]
&\leq&\dis K_1|x|\big(g_k(|x|+R)-g_k(|x|)\big)=K_1kRg_k(|x|)\leq\la_kL_k(x),
\ec
where $\la_k:=kK_1R$ with $K_1$ as in (\ref{K01}), and in the last equality we have used (\ref{diftrick}).
\epro

\bpro[of Theorem~\ref{T:poly}]
The fact that $G$ maps $\Fi_\infty$ into itself follows from Lemma~\ref{L:GFk}. By formula (\ref{Plat}) in the appendix,
\be
P_tL_k\leq e^{\la_kt}L_k\qquad(t\geq 0, \ k \geq 0),
\ee
which implies that $P_t$ maps $\Fi_k$ into itself for each $t\geq 0$ and $k\geq 0$. To prove the remaining statements we will apply Theorem~\ref{T:Lbp} from the appendix. Since $f\in\Fi_\infty$, there exists a $k\geq 0$ such that $f\in\Fi_k$. Then $Gf\in\Fi_{k+1}$ by Lemma~\ref{L:GFk}. It follows that both $f$ and $Gf$ are $L_{k+1}$-bounded in the terminology of Theorem~\ref{T:Lbp}. This allows us to apply that theorem with $L=L_{k+1}$ and $L'=L_{k+2}$, provided we show that
\be\label{diabd}
\sup\big\{\!-G(x,x):L_{k+2}(x)\leq NL_{k+1}(x)\big\}<\infty\quad\forall N<\infty.
\ee 
Here
\be
\frac{L_{k+2}(x)}{L_{k+1}(x)}=\frac{g_{k+2}(|x|)+1}{g_{k+1}(|x|)+1}
=\frac{1+\prod_{i=0}^{k+1}(|x|+iR)}{1+\prod_{i=0}^k(|x|+iR)}
\ee
tends to $\infty$ as $|x|\to\infty$, so for a given $N<\infty$ one has $L_{k+2}(x)\leq NL_{k+1}(x)$ only for finitely many values of $|x|$. In view of this, Lemma~\ref{L:Gdia} implies that (\ref{diabd}) is satisfied and the proof is complete.
\epro

\subsection{Monotonicity}

We continue to work in the set-up of Subsection~\ref{S:add}, so $\Si:=\Lt^\La$ where $\Lt$ is a finite lattice and $\La$ a countable set, and we are interested in interacting particle systems $(X_t)_{t\geq 0}$ with generator of the form (\ref{Gf}) where $\Gi$ is a set whose elements are additive local maps $m\cn\Si\to\Si$ and~$(r_m)_{m\in\Gi}$ are nonnegative rates satisfying (\ref{sumad}). In the present subsection, we won't need~(\ref{Rbd}).

It is easy to see that each additive map $m$ is monotone in the sense that $x\leq y$ implies~$m(x)\leq m(y)$. It is well-known that this implies that the process $(X_t)_{t\geq 0}$ has an invariant law $\ov\nu_X$ such that
\be
\P^{\un\top}\big[X_t\in\,\cdot\,\big]\Asto{t}\ov\nu_X,
\ee
where $\Rightarrow$ denotes weak convergence of probability measures on $\Si$. For a proof in the case that $\Lt=\{0,1\}$ see \cite[Thm~III.2.3]{Lig85} or \cite[Thm~5.13]{Swa26IPS}. The proof for general finite lattices is the same. The invariant law $\ov\nu_X$ is called the \emph{upper invariant law} of $(X_t)_{t\geq 0}$. Similarly, $(Y_t)_{t\geq 0}$ has an upper invariant law $\ov\nu_Y$. Generalising the definitions in (\ref{tetdef}) and (\ref{vtetdef}), we set
\be
\tet^a_i:=\P^{e^a_i}\big[X_t\neq\un 0\ \forall t\geq 0\big]
\quand
\vtet^a_i:=\int\ov\nu(\di x)1_{\txt\{x(i)\not\leq a\big\}}\qquad(a\in\Lt,\ i\in\La).
\ee
We let $\hat\tet^b_i$ and $\hat\vtet^b_i$ denote the same quantities for the dual process $(b\in\Lt',\ i\in\La)$. We claim that
\be\label{tetvar}
\vtet^a_i=\hat\tet^{a'}_i\quand\hat\vtet^b_i=\tet^{b'}_i\qquad(a\in\Lt,\ b\in\Lt',\ i\in\La).
\ee
To see this, define $f^a_i\in\Si$ by $f^a_i(j):=a$ if $i=j$ and $:=\top$ otherwise. Then $(f^a_i)'=e^{a'}_i$ and duality (\ref{dual}) implies
\be
1-\vtet^a_i
=\lim_{t\to\infty}\P^{\un\top}\big[X_t\leq f^a_i\big]
=\lim_{t\to\infty}\P^{\un\top}\big[e^{a'}_i\leq X'_t\big]
=\lim_{t\to\infty}\P^{e^{a'}_i}\big[Y_t\leq\un 0\big]
=1-\hat\tet^{a'}_i.
\ee
The proof that $\hat\vtet^b_i=\tet^{b'}_i$ is the same. Formula (\ref{fininf}) is a special case of (\ref{tetvar}), where in this case $\Lt=\Lt'=\{0,1,2\}$ and $a':=2-a$, and the quantities $\tet^a_i$ and $\vtet^a_i$ do not depend on $i$ by translation invariance. For two-stage contact processes, we have the following lemma, which has already been announced in the introduction.

\bl[Equivalent definitions of survival]
Fix a function $\la\cn \La \times \La\to [0,\infty)$ as in Section \ref{S:2stage}. For\label{L:tet12} $\ga>0$ and $\de_1,\de'_1,\de_2\geq 0$ one has $\tet^1(\la, \ga,\de_1,\de'_1,\de_2)>0$ if and only if $\tet^2(\la, \ga,\de_1,\de'_1,\de_2)>0$ and $\vtet^0(\la,\ga,\de_1,\de'_1,\de_2)>0$ if and only if $\vtet^1(\la, \ga,\de_1,\de'_1,\de_2)>0$.
\el

\bpro
Monotonicity of the maps $\Xb_{s,t}$ gives $\tet^1\leq\tet^2$. On the other hand, $\tet^2>0$ implies $\tet^1>0$ by the assumption $\ga>0$. The statement for $\vtet^0$ and $\vtet^1$ now follows from (\ref{fininf}).
\epro

We will not only need monotonicity in the initial state but also monotonicity in the rates. We write $\Gi=\Gi_+\cup\Gi_-\cup\Gi_{+-}$ where
\be\label{Gipm}
\Gi_+:=\big\{m\in\Gi:x\leq m(x)\ \forall x\in\Si\big\},\quad
\Gi_-:=\big\{m\in\Gi:m(x)\leq x\ \forall x\in\Si\big\},
\ee
and $\Gi_{+-}:=\Gi\beh(\Gi_+\cup\Gi_-)$.

\bl[Monotone coupling]
Let\label{L:moncoup} $\Gi$ be a set whose elements are additive local maps $m\cn\Si\to\Si$. Let $(r_m)_{m\in\Gi}$ and $(r'_m)_{m\in\Gi}$ be nonnegative rates satisfying (\ref{sumad}) such that
\be
r_m\leq r'_m\ \forall m\in\Gi_+,\quad
r_m\geq r'_m\ \forall m\in\Gi_-,\quand
r_m=r'_m\ \forall m\in\Gi_{+-}.
\ee
Then it is possible to couple two graphical representations $\om$ and $\om'$ with rates $(r_m)_{m\in\Gi}$ and $(r'_m)_{m\in\Gi}$ in such a way that the corresponding stochastic flows satisfy
\be\label{moncoup}
\Xb_{s,t}(x)\leq\Xb'_{s,t}(x)\qquad(s\leq t,\ x\in\Si). 
\ee
\el

\bpro
Set $r^\circ_m:=r_m\wedge r'_m$, $r^+_m:=(r'_m-r_m)\vee 0$, and $r^-_m:=(r_m-r'_m)\vee 0$ $(m\in\Gi)$ and let $\om^\circ$, $\om^+$, and $\om^-$ be independent Poisson point sets with intensity as in (\ref{rho}) and $r_m$ replaced by $r^\circ_m$, $r^+_m$, and $r^-_m$, respectively. Then we claim that setting $\om':=\om^\circ\cup\om^+$ and $\om:=\om^\circ\cup\om^-$ does the trick. Indeed, these are Poisson point sets with the right rates such that $\om'$ compared to $\om$ has more elements of the form $(m,t)$ with $m\in\Gi_+$ and fewer elements of the form $(m,t)$ with $m\in\Gi_-$. Using the fact that additive maps are monotone, it is intuitively clear that this implies (\ref{moncoup}), at least on finite grids, where $\Xb_{s,t}$ is the concatenation of finitely many local maps. To also prove the statement on infinite grids one can apply an approximation result such as \cite[Prop~4.23]{Swa26IPS}.
\epro

As an immediate consequence of Lemma~\ref{L:moncoup}, for $(\la,\ga,\de_1,\de'_1,\de_2)$-contact processes, we obtain:

\bl[Monotonicity in the rates]\label{lemmamonotonicity}
The functions $\tet^1$ and $\tet^2$ are nondecreasing in $\ga$ and nonincreasing in $\de_1,\de'_1$, and $\de_2$.
\el

\subsection{The exponential growth rate}\label{S:exp}

In this subsection we prove a more general version of Lemma~\ref{L:exprate}, which holds for additive interacting particle systems. As in the previous subsections, $\Si:=\Lt^\La$ where $\Lt$ is a finite lattice and $\La$ a countable set. We are interested in the interacting particle system $(X_t)_{t\geq 0}$ with generator of the form (\ref{Gf}) where $\Gi$ is a set whose elements are additive local maps $m\cn\Si\to\Si$ and $(r_m)_{m\in\Gi}$ are nonnegative rates satisfying (\ref{sumad}). In the present subsection, we moreover assume that $\La$ is a countable group and that the expression (\ref{Gf}) is invariant under the left shift of the group on itself, as we explain now.

For $j\in\La$ we define a shift operator $T_j\cn\Si\to\Si$ by
\be\label{shift}
(T_jx)(i):=x(j^{-1}i)\qquad(i\in\La).
\ee
It is easy to check that $T_iT_j=T_{ij}$,
the inverse of $T_j$ is $T^{-1}_j=T_{j^{-1}}$, and $T_je^a_i=e^a_{ji}$ $(i,j\in\La)$.
If $m\cn\Si\to\Si$ is a map, then for $j\in\La$ we define a shifted map $\Ti_jm\cn\Si\to\Si$ by
\be
\Ti_jm:=T_j\circ m\circ T^{-1}_j.
\ee
For example, for the maps $m^{1\to 2}_i$ and $m^\ast_{i,j}$ defined in (\ref{mapi}) and (\ref{mapij}) one can check that
\be
\Ti_km^{1\to 2}_i=m^{1\to 2}_{ki}
\quand
\Ti_km^\ast_{i,j}=m^\ast_{ki,kj}.
\ee
We assume from now on that the random mapping representation (\ref{Gf}) is translation invariant in the sense that
\be\label{homog}
m\in\Gi\quad\mbox{implies}\quad\Ti_km\in\Gi\quand r_{\Ti_km}=r_m.
\ee
The following formula implies (\ref{exprate}) as a special case. We call $r$ the \emph{exponential growth rate} of the additive interacting particle system $(X_t)_{t\geq 0}$.

\bl[The exponential growth rate]
There\label{L:exp} exists a constant $r\in\R$ such that
\be\label{exp}
r=\lim_{t\to\infty}\frac{1}{t}\log\E^{e^\top_0}\big[|X_t|]=\inf_{t>0}\frac{1}{t}\log\E^{e^\top_0}\big[|X_t|].
\ee
\el

\bpro
Set $\pi_t:=\E^{e^\top_0}\big[|X_t|]$ which is finite by (\ref{expbd}). Using translation invariance as well as the fact that the stochastic flow $(\Xb_{s,t})_{s\leq t}$ is monotone and stationary with independent increments, we see that
\bc
\dis\pi_{s+t}&=&\dis\sum_{j\in\La}\P\big[\Xb_{0,s+t}(e^\top_0)(j)\neq 0\big]\\[5pt]
&\leq&\dis\sum_{j\in\La}\P\big[\exists i\in\La\mbox{ s.t.\ }\Xb_{0,s}(e^\top_0)(i)\neq 0\quand\Xb_{s,s+t}(e^\top_i)(j)\neq 0\big]\\[5pt]
&\leq&\dis\sum_{i,j\in\La}\P\big[\Xb_{0,s}(e^\top_0)(i)\neq 0\big]\cdot\P\big[\Xb_{s,s+t}(e^\top_i)(j)\neq 0\big]\\[5pt]
&=&\dis\sum_{i,j\in\La}\P\big[\Xb_{0,s}(e^\top_0)(i)\neq 0\big]\cdot\P\big[\Xb_{0,t}(e^\top_0)(i^{-1}j)\neq 0\big]=\pi_s\pi_t\quad(s,t\geq 0).
\ec
This shows that $t\mapsto\pi_t$ is submultiplicative and hence $t\mapsto\log\pi_t$ is subadditive. Fekete's lemma \cite{Fek23} now implies that $\lim_{t\to\infty}\frac{1}{t}\log\pi_t=\inf_{t>0}\frac{1}{t}\log\pi_t$. A priori, this may be $-\infty$, but using (\ref{sumad}) it is easy to see that $\P[\Xb_{0,t}(e^\top_0)(0)=\top]\geq e^{-Kt}$ for some $K<\infty$, which implies $t^{-1}\log\pi_t\geq-K$ for all $t>0$.
\epro

\bl[The dual exponential growth rate]
One\label{L:duexp} has $\E^{e^\top_0}\big[|X_t|]=\E^{e^\top_0}\big[|Y_t|]$ $(t\geq 0)$ so the exponential growth rates of $(X_t)_{t\geq 0}$ and $(Y_t)_{t\geq 0}$ are the same.
\el

\bpro
This follows from duality (\ref{dual}) and translation invariance by writing
\be\ba{l}
\dis\E^{e^\top_0}\big[|Y_t|\big]
=\sum_{i\in\La}\P^{e^\top_0}\big[Y_t\not\leq(e^\top_i)'\big]
=\sum_{i\in\La}\P^{e^\top_i}\big[X_t\not\leq(e^\top_0)'\big]\\[5pt]
\dis\quad=\sum_{i\in\La}\P^{e^\top_0}\big[X_t\not\leq(e^\top_{i^{-1}})'\big]
=\E^{e^\top_0}\big[|X_t|\big].
\ec
\epro

\bpro[of Lemma~\ref{L:exprate}]
Immediate from Lemmas \ref{L:exp} and \ref{L:duexp}.
\epro

For later use, we prove the following simple consequence of Lemma~\ref{L:exp}.

\bl[Exponential bound]
For\label{L:exp2} each $x\in\Si_{\rm fin}$ and $\eps>0$, there exists a constant $C<\infty$ such that
\be\label{exp2}
\E^x\big[|X_t|\big]\leq C\ex{(r+\eps)t}\qquad(t\geq 0).
\ee
\el

\bpro
By additivity of $\Xb_{0,t}$ and translation invariance we have $\E^x[|X_t|]\leq|x|\E^{e^\top_0}[|X_t|]$, so Lemma~\ref{L:exp} implies that for each $\eps>0$ one has $\E^x[|X_t|]\leq e^{(r+\eps)t}$ for all $t$ large enough. Since $t\mapsto\E^x[|X_t|]$ is bounded on bounded time intervals (see (\ref{expbd})), it follows that (\ref{exp2}) holds for $C$ large enough.
\epro

\section{Eigenmeasures and sharpness of the phase transition}

\subsection{Locally finite measures}

In this section we prove Theorem~\ref{T:sharp}. A key tool in the proof will be \emph{eigenmeasures}, which are locally finite measures on the space of nonzero configurations that under the forward equation of an additive particle system are mapped into a constant multiple of themselves. For contact processes on groups, the existence of such eigenmeasures was first established in \cite{Swa09}. We will need to generalise some results that were proved in \cite{Swa09} for $\Lt=\{0,1\}$ to general finite lattices.

As in the previous sections, $\Si:=\Lt^\La$ where $\Lt$ is a finite lattice and $\La$ a countable set. Since $\Si$ is compact, $\Si_+:=\Si\beh\{\un 0\}$ is a locally compact space. A measure $\mu$ on $\Si_+$ is \emph{locally finite} if $\mu(C)<\infty$ for all compact $C\sub\Si_+$. We equip $\Si_+$ with the Borel-\si-field. We start with a pretty straightforward generalisation of \cite[Lemma~3.1]{Swa09}.

\bl[Locally finite measures]
For\label{L:locfin} a measure $\mu$ on $\Si_+$, the following are equivalent: (i) $\mu$ is locally finite, (ii) $\int\mu(\di x)1_{\{x(i)\neq 0\}}<\infty$ for all $i\in\La$, (iii) $\int\mu(\di x)\psi(x,y)<\infty$ for all~$y\in\Si'_{\rm fin}$.
\el

\bpro
For each $y\in\Si'_{\rm fin}$ the set $C_y:=\{x\in\Si:x\not\leq y'\}$ is a closed subset of $\Si$ and therefore compact. Moreover $C_y\sub\Si_+$ so if $\mu$ is locally finite then $\int\mu(\di x)1_{\{x\not\leq y'\}}=\mu(C_y)<\infty$. In particular, setting $y=e^\top_i$ we see that $\int\mu(\di x)1_{\{x(i)\neq 0\}}<\infty$. This proves (i)$\volgt$(iii)$\volgt$(ii). To complete the proof we will show that (ii)$\volgt$(i). Let $C\sub\Si_+$ be compact. Let $\De_n$ be finite sets increasing to $\La$ and set $C_n:=\{x\in\Si_+:x(i)=0 \ \forall i\in\De_n\}$. The sets $C_n$ are closed and decrease to the empty set so we can choose $n$ large enough such that $C_n\cap C=\emptyset$. Then $\mu(C)\leq\sum_{i\in\De_n}\int\mu(\di x)1_{\{x(i)\neq 0\}}<\infty$ proving that $\mu$ is locally finite.
\epro

We will also need the following simple lemma.

\bl[Measure-determining functions]
If\label{L:measdet} two locally finite measures $\mu$ and $\nu$ on $\Si_+$ satisfy $\int\mu(\di x)\psi(x,y)=\int\nu(\di x)\psi(x,y)$ for all $y\in\Si'_{\rm fin}$, then $\mu=\nu$.
\el

\bpro
Since $\phi(x,y_1)\phi(x,y_2)=1_{\{y_1\leq x'\}}1_{\{y_2\leq x'\}}=1_{\{y_1\vee y_2\leq x'\}}$, the linear span of the functions $\phi(\,\cdot\,,y)$ with $y\in\Si'_{\rm fin}$ forms an algebra that separates points, so by the Stone--Weierstrass theorem these functions are dense in the space of all continuous functions on $\Si$, equipped with the supremum-norm. It follows that if two finite measures $\mu'$ and $\nu'$ on $\Si$ satisfy $\int\mu'(\di x)\phi(x,y)=\int\nu'(\di x)\phi(x,y)$ for all $y\in\Si'_{\rm fin}$, then $\mu'=\nu'$.

Now let $\mu$ and $\nu$ be as in the lemma. Let $\De_m$ be finite sets increasing to $\La$, and define $y_m\in\Si'_{\rm fin}$ by $y_m(i):=\top$ if $i\in\De_m$ and $:=0$ otherwise. Then the $\psi(\,\cdot\,,y_m)$ are nonnegative, continuous, compactly supported functions on $\Si_+$. Since $\psi(x,y_m)=1_{\{\exists i\in\De_m\mbox{ s.t.\ }x(i)\neq 0\}}$, the functions $\psi(\,\cdot\,,y_m)$ increase to $1$ as $m\to\infty$. Since for each $y\in\Si'_{\rm fin}$
\be
\int\!\mu(\di x)\psi(x,y_m)\phi(x,y)=\int\!\mu(\di x)\psi(x,y_m\vee y)-\int\!\mu(\di x)\psi(x,y)=\int\!\nu(\di x)\psi(x,y_m)\phi(x,y),
\ee
the weighted measures $\mu(\di x)\psi(x,y_m)$ and $\nu(\di x)\psi(x,y_m)$ are equal for each $m$, so letting $m\to\infty$ we see that $\mu=\nu$.
\epro

We equip the space of locally finite measures on $\Si_+$ with the vague topology, in which a sequence of locally finite measures $\mu_n$ converges to a limit $\mu$, denoted $\mu_n\To\mu$, if $\int\mu_n(\di x)f(x)\to\int\mu(\di x)f(x)$ for each continuous compactly supported real function $f$. It is well-known that the topology of vague convergence is metrisable. The next lemma is a pretty straightforward generalisation of \cite[Lemma~3.2]{Swa09}.

\bl[Vague convergence]
Let\label{L:vague} $\mu_n,\mu$ be locally finite measures on $\Si_+$. Then the $\mu_n$ converge vaguely to $\mu$ if and only if $\int\mu_n(\di x)\psi(x,y)\to\int\mu(\di x)\psi(x,y)$ for all $y\in\Si'_{\rm fin}$. The sequence $\mu_n$ is relatively compact in the topology of vague convergence if and only if $\sup_n\int\mu_n(\di x)\psi(x,y)<\infty$ for each $y\in\Si'_{\rm fin}$.
\el

\bpro
For each $y\in\Si'_{\rm fin}$, the function $\Si_+\ni x\mapsto\psi(x,y)$ is continuous with compact support, so the conditions for convergence and relative compactness are clearly necessary. To see that they are also sufficient, we start by making a few general observations.

Let $\Xc$ be a locally compact, metrisable space, and let $f_m$ be nonnegative, continuous, compactly supported functions on $\Xc$ whose supports increase to $\Xc$ as $m\to\infty$. Then a sequence of locally finite measures $\mu_n$ on $\Xc$ converges vaguely to a limit if and only if the weighted measures $f_m(x)\mu_n(\di x)$ converge weakly to $f_m(x)\mu(\di x)$ for each $m$. By Prokhorov's theorem, a sequence of finite measures on a compact space is relatively compact in the topology of weak convergence if and only if their total measures are bounded.

Applying this to the functions $f_m:=\psi(\,\cdot\,,y_m)$ with $y_m$ as in the proof of Lemma~\ref{L:measdet}, we see that if $\sup_n\int\mu_n(\di x)\psi(x,y_m)<\infty$ for each $m$, then the sequence $\mu_n$ is relatively compact in the topology of vague convergence. If moreover $\int\mu_n(\di x)\psi(x,y)\to\int\mu(\di x)\psi(x,y)$ for all $y\in\Si_{\rm fin}$, then by Lemma~\ref{L:measdet} each vague cluster point of the $\mu_n$ is equal to $\mu$, which implies $\mu_n\To\mu$.
\epro

\subsection{Evolution of locally finite measures}

In the present subsection, we will show that it is possible to start additive particle systems in a locally finite measure on $\Si_+$. This generalises theory developed in \cite{Swa09} for the contact process.

We will use the following notation. If $M$ is a measurable space, $\mu$ is a measure on $M$, $f\cn M\to[0,\infty]$ is a measurable function, and $K,L$ are measurable kernels on $M$, then we write
\be
\mu K:=\int\mu(\di x)K(x,\,\cdot\,),\quad
Kf(x):=\int K(x,\di y)f(y),\quad
(KL)(x,\,\cdot\,):=\int K(x,\di y)L(y,\,\cdot).
\ee
We work in the general set-up of Subsection~\ref{S:exp}, so $(X_t)_{t\geq 0}$ is a Markov process with state space $\Si$ and generator of the form (\ref{Gf}) where the rates $(r_m)_{m\in\Gi}$ satisfy the summability condition (\ref{sumad}) and are translation invariant in the sense of (\ref{homog}). We define measurable kernels $P^+_t$ $(t\geq 0)$ on $\Si_+$ by
\be\label{Pt}
P^+_t(x,\Ai):=\P^x\big[X_t\in\Ai\big]\qquad\forall t\geq 0,\ x\in\Si_+,\ \mbox{measurable}\ \Ai\sub\Si_+.
\ee
Note that $P^+_t$ is the transition kernel of $(X_t)_{t\geq 0}$, restricted to the smaller space $\Si_+$, so $P^+_t(x,\,\cdot\,)$ is a subprobability measure but in general not a probability measure. Since $\un 0$ is a trap for the process $(X_t)_{t\geq 0}$, it is easy to see that the kernels $(P^+_t)_{t\geq 0}$ form a semigroup in the sense that $P^+_sP^+_t=P^+_{s+t}$ $(s,t\geq 0)$. We define $(\hat P^+_t)_{t\geq 0}$ as in (\ref{Pt}) but with $(X_t)_{t\geq 0}$ replaced by the dual process $(Y_t)_{t\geq 0}$. 

We say that a measure $\mu$ on $\Si$ or on $\Si_+$ is \emph{homogeneous} if
\be
\mu\circ T_i=\mu\qquad(i\in\La),
\ee
where $T_i$ is the shift operator defined in (\ref{shift}). The following lemma generalises \cite[Lemma~3.3]{Swa09}. Below $\To$ denotes vague convergence on $\Si_+$.

\bl[Evolution of locally finite measures]
If\label{L:finevol} $\mu$ is a homogeneous, locally finite measure on $\Si_+$, then so is $\mu P^+_t$ for each $t\geq 0$. If $\mu_n,\mu$ are homogeneous, locally finite measures on $\Si_+$ such that $\mu_n\To\mu$, then $\mu_nP^+_t\To\mu P^+_t$ for all $t\geq 0$.
\el

\bpro
The local finiteness of $\mu P^+_t$ follows with the help of Lemma~\ref{L:locfin} by writing, for arbitrary~$i\in\La$,
\be\ba{l}\label{Pdu}
\dis\int\!(\mu P^+_t)(\di x)1_{\{x(i)\neq 0\}}
=\int\!(\mu P^+_t)(\di x)1_{\{x\not\leq(e^\top_i)'\}}
=\int\!\mu(\di x)\E^x\big[\psi(X_t,e^\top_i)\big]\\[5pt]
\dis\quad=\int\!\mu(\di x)\E^{e^\top_i}\big[\psi(x,Y_t)\big]
=\int\!\hat P^+_t(e^\top_i,\di z)\int\!\mu(\di x)\psi(x,z)\\[5pt]
\dis\quad\leq\int\!\hat P^+_t(e^\top_i,\di z)\!\sum_{i:\,z(i)\neq 0}\int\!\mu(\di x)1_{\{x(i)\neq 0\}}
=\E^{e^\top_i}\big[|Y_t|\big]\int\!\mu(\di x)1_{\{x(0)\neq 0\}}<\infty,
\ec
where we have used duality and the homogeneity of $\mu$ and the final expression is finite by (\ref{expbd}) and the local finiteness of $\mu$. The homogeneity of $\mu P^+_t$ follows from the homogeneity of $\mu$ and the translation invariance of the rates. To show that $\mu_n\To\mu$ implies $\mu_nP^+_t\To\mu P^+_t$, by Lemma~\ref{L:vague}, it suffices to show that
\be
\int\!(\mu_n P^+_t)(\di x)\psi(x,y)\to\int\!(\mu P^+_t)(\di x)\psi(x,y)\qquad(y\in\Si'_{\rm fin}).
\ee
By our previous calculation
\be\label{munP}
\int\!(\mu_n P^+_t)(\di x)\psi(x,y)=\int\!\hat P^+_t(y,\di z)\int\!\mu_n(\di x)\psi(x,z).
\ee
Here $\mu_n(\di x)\psi(x,z)$ converges pointwise for each $z\in\Si'_{\rm fin}$ by our assumption $\mu_n\To\mu$ and Lemma~\ref{L:vague}. This function can moreover be bounded from above by
\be
|z|\int\!\mu_n(\di x)1_{\{x(0)\neq 0\}}
\ee
where the integral is bounded uniformly in $n$ since $\mu_n\To\mu$. By dominated convergence, using (\ref{expbd}), we see that the integral in (\ref{munP}) converges to the same thing with $\mu_n$ replaced by $\mu$.
\epro

\subsection{Eigenmeasures}

We continue to work in the general set-up of Subsection~\ref{S:exp}, so $(X_t)_{t\geq 0}$ is an additive interacting particle system with state space $\Si$ and generator of the form (\ref{Gf}) where the rates $(r_m)_{m\in\Gi}$ satisfy the summability condition (\ref{sumad}) and are translation invariant in the sense of (\ref{homog}). Following \cite{Swa09}, we define an \emph{eigenmeasure} of the process $(X_t)_{t\geq 0}$ to be a nonzero, locally finite measure $\mu$ on $\Si_+$ for which there exists a constant $r\in\R$ such that
\be\label{eig}
\mu P^+_t=e^{rt}\mu\qquad(t\geq 0).
\ee
We call $r$ the associated \emph{eigenvalue}. We will prove the following proposition. For contact processes on groups, this has been proved before in \cite[Prop.~1.4]{Swa09}. The proof below closely follows the proof there, with some modifications due to the more complicated definition of the exponential growth rate.

\bp[Existence of a homogeneous eigenmeasure]
The\label{P:eigmeas} process $(X_t)_{t\geq 0}$ has a homogeneous eigenmeasure whose eigenvalue $r$ is the exponential growth rate defined in Lemma~\ref{L:exp}.
\ep

\bpro[of Proposition~\ref{P:eigmeas}]
Let $\de_{e^\top_i}$ denote the delta-measure on the configuration $e^\top_i$ and let
\be
\mu:=\sum_{i\in\La}\de_{e^\top_i}.
\ee
 Set $\mu_t:=\mu P_t^+$ ($t\geq 0$). By Lemma~\ref{L:locfin}, $\mu$ is a locally finite measure on $\Si_+$, so by Lemma~\ref{L:finevol} the same is true for $\mu_t$ for all $t\geq 0$. The calculation in (\ref{Pdu}) gives
\be\ba{l}\label{mtins}
\dis\int\!\mu_t (\di x)1_{\{x(0)\neq 0\}}
=\int\!\hat P^+_t(e^\top_0,\di z)\int\!\mu(\di x)\psi(x,z)
=\int\!\hat P^+_t(e^\top_0,\di z)\sum_{i\in\La}\psi(e^\top_i,z)\\[5pt]
\dis\quad=\int\hat P^+_t(e^\top_0,z)\sum_{i\in\La}1_{\txt\{z\not\leq(e^\top_i)'\}}
=\int\hat P^+_t(e^\top_0,z)|z|=\E^{e^\top_0}\big[|Y_t|\big].
\ec
 Let $\hat\mu_\al$ denote the Laplace transform of $\mu_t$, that is
\be
\hat\mu_\al:=\int_0^\infty\!\mu_t\,e^{-\al t}\di t\qquad(\al\in\R).
\ee
As in the proof of Lemma~\ref{L:exp} we set $\pi_t:=\E^{e^\top_0}[|X_t|]$, which equals $\E^{e^\top_0}[|Y_t|]$ by Lemma~\ref{L:duexp}. Lemma~\ref{L:exp} implies that for each $\eps>0$ one has $\pi_t\leq e^{(r+\eps)t}$ for all $t$ large enough, so (\ref{mtins}) implies
\be
\hat\pi_\al:=\int\!\hat\mu_\al(\di x)1_{\{x(0)\neq 0\}}=\int_0^\infty\!\pi_t\,e^{-\al t}\di t<\infty\quad(\al>r).
\ee
By spatial homogeneity and Lemma~\ref{L:locfin} we conclude that $\hat\mu_\al$ is a locally finite measure for all $\al>r$. Lemma~\ref{L:exp} moreover tells us that $t^{-1}\log\pi_t\geq r$ for all $t>0$ so by monotone convergence
\be\label{tinf}
\lim_{\al\down r}\hat\pi_\al=\int_0^\infty\!\pi_t\,e^{-rt}\di t\geq\int_0^\infty\!e^{rt}\,e^{-rt}\di t=\infty.
\ee
Set $\ov\mu_\al:=\hat\pi_\al^{-1}\hat \mu_\al$ $(\al>r)$. These measures are homogeneous, locally finite, and normalised such that $\int\ov\mu_\al(\di x)1_{\{x(0)\neq 0\}}=1$. Similarly to what we did in (\ref{Pdu}), for general $y\in\Si'_{\rm fin}$ we can estimate
\be
\int\ov\mu_\al(\di x)\psi(x,y)\leq\sum_{i:\,y(i)\neq 0}\int\ov\mu_\al(\di x)1_{\{x(i)\neq 0\}}=|y|,
\ee
so we can apply Lemma~\ref{L:vague} to conclude that if $\al_n>r$ satisfy $\al_n\to r$, then by going to a subsequence we can assume that $\ov\mu_{\al_n}\To\ov\mu_r$ for some homogeneous, locally finite $\ov\mu_r$. Obviously $\int\ov\mu_r(\di x)1_{\{x(0)\neq 0\}}=1$ so $\ov\mu_r$ is nonzero. Moreover, for each $t\geq 0$,
\bc
\dis\ov\mu_rP^+_t
&=&\dis\lim_{n\to\infty}\hat\pi_{\al_n}^{-1}\int_0^\infty\! e^{-\al_ns}\di s\int\mu(\di x)\int\!P^+_s(x,\di z)\,P^+_t(z,\,\cdot\,)\\[5pt]
&=&\dis e^{rt}\lim_{n\to\infty}\hat\pi_{\al_n}^{-1}\int_0^\infty\! e^{-\al_n(s+t)}\di s\int\mu(\di x)P^+_{s+t}(x,\,\cdot\,)\\[5pt]
&=&\dis e^{rt}\lim_{n\to\infty}\hat\pi_{\al_n}^{-1}\Big(\int_0^\infty\! e^{-\al_ns}\di s\int\mu(\di x)P^+_s(x,\,\cdot\,)-\int_0^t\! e^{-\al_ns}\di s\int\mu(\di x)P^+_s(x,\,\cdot\,)\Big)\\[5pt]
&=&\dis e^{rt}\ov\mu_r-e^{rt}\lim_{n\to\infty}\hat\pi_{\al_n}^{-1}\int_0^t\! e^{-\al_ns}\di s\int\mu(\di x)P^+_s(x,\,\cdot\,).
\ec
Using (\ref{mtins}) and, crucially, (\ref{tinf}), we see that the second term converges vaguely to the zero measure, so we conclude that $\ov\mu_r$ is a homogeneous eigenmeasure with eigenvalue $r$.
\epro

\subsection{Eigenfunctions}

We continue to work in the general set-up of Subsection~\ref{S:exp}, so $(X_t)_{t\geq 0}$ is an additive interacting particle system with state space $\Si$ and generator of the form (\ref{Gf}) where the rates $(r_m)_{m\in\Gi}$ satisfy the summability condition (\ref{sumad}) and are translation invariant in the sense of (\ref{homog}). We let $(Y_t)_{t\geq 0}$ denote the dual process defined in Subsection~\ref{S:add}.

\bl[Eigenfunctions of the generator]
Let\label{L:lemmah} $\mu$ be a homogeneous eigenmeasure of the dual process $(Y_t)_{t\geq 0}$ and let $h\cn\Si_{\rm fin}\to\half$ be defined as
\be
h(x):=\int_{\Si_+}\!\mu(\di y)\,\psi(x,y)\qquad(x\in\Si_{\rm fin}).
\ee
Then, using notation as in (\ref{ddef}) and (\ref{shift}),
\begin{enumerate}
\item $h\circ T_i=h\quad(i\in\La)$,
\item $h(x)\leq h(z)$ for all $x,z\in\Si_{\rm fin}$ with $x\leq z$,
\item $h(x\vee z)\leq h(x)+h(z)$ for all $x,z\in\Si_{\rm fin}$,
\item $h(\un 0)=0$ and $h(e^\top_0)>0$,
\item $\big|h(x)-h(y)\big|\leq Cd(x,y)$ for all $x\in\Si_{\rm fin}$, with $C:=h(e^\top_0)$.
\end{enumerate}
One has
\be \label{eqh}
\E^x\big[h(X_t)\big]=e^{rt}h(x)\qquad(t\geq 0,\ x\in\Si_{\rm fin}).
\ee
If (\ref{Rbd}) holds, then moreover $Gh=rh$, where $G$ is the generator of the process $(X_t)_{t\geq 0}$. 
\el

\bpro
The function $h$ is shift-invariant and monotone. Property~$\text{(iii)}$ follows from the inequality $1_{\{x \vee z\not \leq y'\}}\leq 1_{\{x \not \leq y'\}}+ 1_{\{ z\not \leq y'\}}$ for all $x, z \in \Si, y \in \Si'$. Property~(iv) follows from the observations that $\psi(\un 0,y)=0$ for all $y\in\Si$ while $\int\mu(\di y)\psi(e^\top_0,y)=\int\mu(\di y)1_{\{y(0)\neq 0\}}>0$ by translation invariance and the fact that $\mu$ is nonzero. To prove property~(v), we can without loss of generality assume that there exists a unique $i\in\La$ with $x(i)\neq y(i)$, as the general case then follows by induction. By (ii), $h(v)\leq h(x),h(y)\leq h(w)$ where $v, w \in \mathcal{S}_{\text{fin}}$ are defined by $v(j)=w(j):=x(j)=y(j)$ for $j\neq i$ and $v(i):=0$, $w(i):=\top$. Now (ii) and (iii) give $h(v)\leq h(w)\leq h(v)+h(e^\top_i)$, so using also (i), the claim follows.

Using duality and the fact that $\mu$ is an eigenmeasure of $(Y_t)_{t\geq 0}$, we see that for any $t\geq 0$ and $x\in\Si_{\rm fin}$,
    \begin{equation} 
    \label{eq:eq}
     \begin{split}
        \mathbb E^x[h(X_t)]&= \int \mu(dy)\mathbb E^x[\psi(X_t,y)]= \int \mu(dy)\mathbb E^y[\psi(x, Y_t)]
        \\&=e^{r t} \int \mu(dy) \psi(x,y)=e^{r t} h(x),
    \end{split}
    \end{equation}
proving \eqref{eqh}. Let $(P_t)_{t\geq 0}$ be the semigroup of the process $(X_t)_{t\geq 0}$. Then \eqref{eqh} says that $P_t h= e^{rt}h$. Property~(v) shows that $h$ is of at most linear growth. Together with the assumption (\ref{Rbd}), this allows us to apply Theorem~\ref{T:poly} to conclude that
    \begin{equation*}
        G h =\pbplim_{t\down 0}\frac{1}{t} (P_t h-h)=\pbplim_{t\down 0}\frac{1}{t}(e^{rt}h-h)=rh.
    \end{equation*}
\epro

\subsection{Exponentially growing processes}

In this section we prove the following result, that generalises \cite[Thm~1.2~(d)]{Swa09}, where this is proved for contact processes.

\bp[Exponentially growing processes]
Let\label{P:posr} $(X_t)_{t\geq 0}$ be an additive interacting particle system on a countable group $\La$ whose rates satisfy (\ref{sumad}), (\ref{Rbd}), and (\ref{homog}). Assume that the exponential growth rate $r$ defined in (\ref{exp}) satisfies $r>0$. Then
\be\label{posr}
\P^{e^\top_0}\big[X_t\neq\un 0\ \forall t\geq 0\big]>0.
\ee
\ep

\begin{Proof}
Let $\mu$ be a homogeneous eigenmeasure of the dual process $(Y_t)_{t\geq 0}$ with eigenvalue $r$ and let the function $h$ be defined in terms of $\mu$ as in Lemma \ref{L:lemmah}. Then $h(e^\top_0)>0$ by part~(iv) of that lemma. We will prove (\ref{posr}) more generally with $e^\top_0$ replaced by any $0\neq z\in\Si_{\rm fin}$ such that $h(z)>0$. To do this, it suffices to show that the limit 
\begin{equation*}
   W:= \lim_{t \to \infty} e^{-rt} h(X_t) 
\end{equation*}
exists and satisfies $\mathbb E^z[W ]=h(z)$. Indeed, by Lemma \ref{L:lemmah} (iv), this then implies that
\begin{equation*}
   0< \mathbb P^z [W>0]\leq \mathbb P^z[X_t \neq 0 \ \forall t \geq 0].
\end{equation*}

Let $(\mathcal{F}_t)_{t\ge 0}$ be the  filtration generated by $(X_t)_{t\ge 0}$. By the Markov property of $(X_t)_{t \geq 0}$ and \eqref{eqh}, for any $t\geq s\geq 0$,
\begin{equation*}
    \mathbb E^z[h(X_{t+s})|\mathcal{F}_s]=\mathbb E^{X_s}[h(X_t)]=e^{rt}h(X_s),
\end{equation*}
which implies that the process $(e^{-rt} h(X_t))_{t\geq 0}$ is a martingale with respect to $(\mathcal F_t)_{t\geq 0}$. Every nonnegative martingale converges, so the limit $W=\lim_{t\to \infty} e^{-rt} h(X_t)$ exists a.s. To prove that $\mathbb E^z[W]=h(z)$, it suffices to show that the random variables $(e^{-rt} h(X_t ))_{t\geq 0}$ are uniformly integrable.

By Lemma~\ref{L:lemmah}~(v), the function $h$ is of linear growth, while formula \eqref{eqh} of the same lemma says that $P_t h=e^{rt} h$ $(t\geq 0)$, where $(P_t)_{t\geq 0}$ is the semigroup of the process $(X_t)_{t\geq 0}$. We can then apply Lemma \ref{L:cov} to conclude that the variance of $h(X_t)$ is given by
\begin{equation}\label{eqvar}
    \text{Var}^z(h(X_t))=\int_0^t\mathbb E^z[\Gamma(e^{rs}h, e^{rs} h)(X_{t-s})]ds= \int_0^t \mathbb E^z[\Gamma(h, h)(X_{t-s})] e^{2rs}ds,
\end{equation}
where by \eqref{eqgamma},
\begin{equation}
    \Gamma(h,h)(x)=\sum_{m\in \mathcal{G}} r_m\big(h(m(x))-h(x)\big)^2\qquad(x \in \mathcal{S}_{\rm fin}).
\end{equation}
By (\ref{sumad}) (i) and (iii),
\be
K:=\sup_{i\in\La}\sum_{m\in\Gi}r_m\sum_{j\in\La}1_{\txt\{m[i,j]\neq{\rm Id}[i,j]\}}<\infty.
\ee
Using this, Lemma~\ref{L:lemmah}~(v), and (\ref{Rbd2}), we can estimate, for some $C<\infty$,
\be\ba{l}
\dis\Ga(h,h)(x)\leq C\sum_{m\in\Gi}r_md\big(m(x),x\big)^2
\leq CR^2\sum_{m\in\Gi}r_m1_{\{m(x)\neq x\}}\\[5pt]
\dis\quad\leq CR^2\sum_{m\in\Gi}r_m1_{\{\exists i,j\in\La\mbox{ s.t.\ }x(i)\neq 0,\ m[i,j]\neq{\rm Id}[i,j]\}}\\[5pt]
\dis\quad\leq CR^2\sum_{i\in\La}1_{\{x(i)\neq 0\}}\sum_{m\in\Gi}r_m\sum_{j\in\La}1_{\{m[i,j]\neq{\rm Id}[i,j]\}}
\leq CR^2K|x|.
\ec
Inserting this into \eqref{eqvar} we obtain
\be\label{Varuni}
\text{Var}^z(e^{-rt}h(X_t))
\leq C'e^{-2rt} \int_0^t\mathbb E^z[|X_{t-s}|]e^{2rs} ds
=C'\int_0^t\mathbb E^z[|X_s|]e^{-2rs} ds,
\ee
where $C':=CR^2K$. By Lemma~\ref{L:exp2}, for each $\eps>0$ we can choose $C''<\infty$ such that $\E^z[|X_{t}|]\leq C''e^{(r+\eps)t}$ $(t\geq 0)$. By assumption $r>0$, so choosing $\eps$ small enough such that $r+\eps-2r<0$, we see from (\ref{Varuni}) that $\text{Var}^z(e^{-rt}h(X_t))$ can be bounded from above uniformly in $t\geq 0$, which proves that the random variables $(e^{-rt} h(X_t))_{t\geq 0}$ are uniformly integrable.
\end{Proof}

\subsection{Sharpness of the phase transition} \label{sec:proofTsharp}

In this subsection we prove Theorem~\ref{T:sharp}. We will first prove a result that holds for a general class of additive interacting particle systems and then specialise to the two-stage contact process to obtain Theorem~\ref{T:sharp}.

Let $(X_t)_{t\geq 0}$ be an additive interacting particle system on a countable group $\La$. In addition to (\ref{sumad}), (\ref{Rbd}), and (\ref{homog}), we will need one more assumption. Recall the definition of the sets $\Gi_\pm\sub\Gi$ in (\ref{Gipm}). We will need the assumption
\be\label{updo}
\Gi=\Gi_-\cup\Gi_+.
\ee
In words, this says that the local maps $m\in\Gi$ come in two types: those whose application can only increase a configuration, and those that can only decrease it. In particular, (\ref{updo}) implies that there are no jumps between incomparable states, a condition known to be equivalent to preservation of positive correlations.\footnote{For finite lattices, the fact that this condition is sufficient and necessary is proved in \cite[Thm~II.2.14]{Lig85} and the discussion below it. For infinite lattices, sufficiency is proved in \cite[Thm~5.11]{Swa26IPS} .} We write
\be
\rb^-:=(r_m)_{m\in\Gi_-}\quand\rb^+:=(r_m)_{m\in\Gi_+}.
\ee
We let $r=r(\rb^-,\rb^+)$ denote the exponential growth rate defined in (\ref{exp}) and we let
\be
\tet=\tet(\rb^-,\rb^+):=\P^{e^\top_0}\big[X_t\neq\un 0\ \forall t\geq 0\big]
\ee
the \emph{survival probability} of the process with rates $\rb^-$ and $\rb^+$. Since multiplying all rates by a constant $a$ has the effect of speeding up time by a factor $a$,
\be\label{timscal}
r(a\rb^-,a\rb^+)=ar(\rb^-,\rb^+)\quand\tet(a\rb^-,a\rb^+)=\tet(\rb^-,\rb^+)\qquad(a>0).
\ee
We will prove the following result.

\bt[Sharpness of the phase transition]
Let\label{T:sharp2} $(X_t)_{t\geq 0}$ be an additive interacting particle system on a countable group $\La$ whose rates satisfy (\ref{sumad}), (\ref{Rbd}), (\ref{homog}), and (\ref{updo}). If $r(\rb^-,\rb^+)=0$, then $\tet\big((1-\eps)\rb^-,\rb^+\big)>0$ for all $0<\eps\leq 1$.
\et

Before proving Theorem~\ref{T:sharp2}, we first show how it implies Theorem~\ref{T:sharp}.\med

\bpro[of Theorem~\ref{T:sharp}]
Clearly $r<0$ implies $\tet=0$, so the inclusions $D_{\rm fin}\sub D_{\rm pos}$ and $\ring D_{\rm fin}\sub\ring D_{\rm pos}$ are trivial. To complete the proof we must show that $\ring D_{\rm pos}\sub\ring D_{\rm fin}$. If $(\ga,\de_1,\de'_1,\de_2)\in\ring D_{\rm pos}$ then by the fact that $\ga>0$ and the fact that $\ring D_{\rm pos}$ is open we can choose $0<\eps\leq\ga/2$ such that $(\ga-2\eps,\de_1+2\eps,\de'_1+2\eps,\de_2+2\eps)\in\ring D_{\rm pos}$. We will show that
\be\label{tetpos}
\tet(\ga-\eps,\de_1+\eps,\de'_1+\eps,\de_2+\eps)>0.
\ee
By monotonicity in the rates (Lemma~\ref{L:moncoup}) this implies that $\tet(\ga-\eps_1,\de_1+\eps_2,\de'_1+\eps_3,\de_2+\eps_4)>0$ for all $\eps_1,\ldots,\eps_4\leq\eps$, from which we conclude that $(\ga,\de_1,\de'_1,\de_2)\in\ring D_{\rm fin}$. Since $(\ga-2\eps,\de_1+2\eps,\de'_1+2\eps,\de_2+2\eps)\in\ring D_{\rm pos}$ we are in one of two cases: I.\ $r(\ga-2\eps,\de_1+2\eps,\de'_1+2\eps,\de_2+2\eps)>0$, and II.\ $r(\ga-2\eps,\de_1+2\eps,\de'_1+2\eps,\de_2+2\eps)=0$. In case~I we apply Proposition~\ref{P:posr} to conclude that $\tet(\ga-2\eps,\de_1+2\eps,\de'_1+2\eps,\de_2+2\eps)>0$ which by monotonicity in the rates (Lemma~\ref{L:moncoup}) implies (\ref{tetpos}). In case~II we apply Theorem~\ref{T:sharp2} which tells us that for each $0<q<1$
\be
\tet\big(\ga-2\eps,q^2(\de_1+2\eps),q^2(\de'_1+2\eps),q^2(\de_2+2\eps)\big)>0.
\ee
By scaling (formula (\ref{timscal})), this implies
\be
\tet\big(q^{-1}(\ga-2\eps),q(\de_1+2\eps),q(\de'_1+2\eps),q(\de_2+2\eps)\big)>0.
\ee
By choosing $q$ sufficiently close to one we can make sure that
\be
q^{-1}(\ga-2\eps)\leq\ga-\eps,\quad
q(\de_1+2\eps)\geq\de_1+\eps,\quad
q(\de'_1+2\eps)\geq\de'_1+\eps,\quad
q(\de_2+2\eps)\geq\de_2+\eps,
\ee
which by monotonicity in the rates (Lemma~\ref{L:moncoup}) implies (\ref{tetpos}).
\epro


\begin{Proof}[of Theorem \ref{T:sharp2}]
Suppose that $r(\rb^-,\rb^+)=0$. By Lemma \ref{L:duexp}, the exponential growth rates of $(X_t)_{t\geq 0}$ and the dual $(Y_t)_{t\geq 0}$ are the same, so Proposition \ref{P:eigmeas} tells us that the dual $(Y_t)_{t\geq 0}$ has a homogeneous eigenmeasure $\mu$ with eigenvalue $0$. We define $h\cn\mathcal{S}_{\text{fin}} \to [0,\infty)$ in terms of $\mu$ as in Lemma~\ref{L:lemmah}. It will be convenient to normalise $\mu$ and $h$ such that
\be\label{normal}
\int_{\Si_+}\!\mu(\di y)\,1_{\{y(0)\neq 0\}}=h(e^\top_0)=1/R,
\ee
where $R$ is the constant from (\ref{Rbd}). The idea of this is that by (\ref{Rbd2}) and Lemma~\ref{L:lemmah}~(v), one then has
\be\label{xmx}
\big|h\big(m(x)\big)-h(x)\big|\leq 1\qquad(x\in\Si_{\rm fin},\ m\in\Gi).
\ee
    
Letting $G$ denote the generator of $(X_t)_{t\geq 0}$, since $\mu$ has eigenvalue zero, Lemma~\ref{L:lemmah} tells us that $Gh = 0$, that is, $h$ is a harmonic function. Let $\varepsilon_1,\varepsilon_2, \varepsilon>0$ be constants to be chosen later. Let $\tilde G$ denote the generator of the process with rates $(1 - \varepsilon_1)\mathbf{r}^-,  (1+\varepsilon_2)\bf{r}^+$ and let
\be
f_{\varepsilon}(x):=\varepsilon^{-1}(1-e^{-\varepsilon h(x)})\qquad(x\in\Si_{\rm fin}).
\ee
We claim that for each $\de>0$, it is possible to choose $\eps_1,\eps_2,\eps>0$ with $\eps_1,\eps_2\leq\de$ such that $\ti Gf_\eps\geq 0$. We first show how this implies the statement of the theorem, and then prove the claim.

Let $(\tilde X_t)_{t\geq 0}$ denote the process with generator $\tilde G$, started in $\tilde X_0=e_0^\top$. We note that by (\ref{poly}), for any function $f$ of polynomial growth, $\ti P_tf=f+\int_0^t\!\di s\,\ti P_s\ti Gf$. In particular, since $\tilde G f_{\varepsilon}\geq 0$ and since $f_\varepsilon$ is a bounded function, this implies that $f_\varepsilon (\tilde X_t)$ is a bounded submartingale that converges to an a.s.\ limit $F_\infty:=\lim_{t\to\infty} f_\varepsilon (\tilde X_t)$ with
\begin{equation*}
    \mathbb E[F_\infty]\geq f_\varepsilon(e_0^\top)=\varepsilon^{-1} (1-e^{-\varepsilon h(e_0^\top)})=\varepsilon^{-1} (1-e^{-\varepsilon/R}),
\end{equation*}
where we used (\ref{normal}). Since $f_\varepsilon(\underline 0)=0$ and $f_\varepsilon\leq\varepsilon^{-1}$, we have
\begin{equation*}
   \mathbb P[\tilde X_t \neq \emptyset \ \forall t\geq 0]\geq \mathbb P [F_\infty>0]\geq \varepsilon \mathbb E[F_\infty]\geq 1-e^{-\varepsilon/R}.
\end{equation*}
By \eqref{timscal}, this implies that
\begin{equation*}
    \theta\left(\frac{(1-\varepsilon_1)}{(1+\varepsilon_2)}\mathbf{r}^-, \mathbf{r}^+\right)=\theta\left((1-\varepsilon_1)\mathbf{r}^-, (1+\varepsilon_2)\mathbf{r}^+ \right)\geq 1-e^{-\varepsilon/R}>0.
\end{equation*}
Since $\eps_1,\eps_2$ can be chosen arbitrarily small, this proves the theorem.

It remains to prove the claim. By \cite[Lemma~3]{Swa18}, one has $\ti Gf_\eps\geq 0$ if and only if
\begin{equation*}
    \ti Gh -\ti H_{\varepsilon} h \geq 0,
\end{equation*}
where
\begin{equation} \label{eqdefH}
    \ti H_{\varepsilon} h(x)=\sum_{y \in \mathcal{S}_{\text{fin}}} \ti G(x,y) \phi_{\varepsilon} (h(y)-h(x)) \quad \text{with} \quad \phi_{\varepsilon}(z):=\varepsilon^{-1} (e^{-\varepsilon z}-1+\varepsilon z),
\end{equation}
and $\ti G(x,y)$ denotes the matrix of $\ti G$, defined as in (\ref{matrix}). Using the fact that $Gh=0$, we see that
 \begin{align*}
    \tilde G h(x)&=\varepsilon_1 \sum_{m \in \mathcal{G}_{-}} r_m\big(h(x)-h(m(x))\big)+\varepsilon_2\sum_{m \in \mathcal{G}_{+}} r_m \big(h(m(x))-h(x)\big).
 \end{align*}
We need to compare this with
 \begin{align*}
     \tilde H_{\varepsilon}h(x)=(1-\varepsilon_1)\sum_{m \in \mathcal{G}_{-}} r_m \phi_\varepsilon\big(h(m(x))-h(x)\big)+(1+\varepsilon_2)\sum_{m \in \mathcal{G}_+} r_m\phi_\varepsilon\big(h(m(x))-h(x)\big).
 \end{align*}
 Since $\phi_\varepsilon (0)=0,\ \phi_\varepsilon' (0)=0$ and $\phi_\varepsilon ''(z)=\varepsilon e^{-\varepsilon z}$, we see that
 \begin{equation*}
     \phi_\varepsilon(z)\leq \ha\varepsilon z^2 \quad (0\leq z \leq 1) \quad \text{and} \quad \phi_\varepsilon(z)\leq \ha\varepsilon e^{\varepsilon}z^2 \quad (-1\leq z\leq 0).
 \end{equation*}
Using (\ref{xmx}) and the fact that $z^2\leq |z|$ for $|z|\leq 1$, we have that
 \begin{align*}
     (\tilde G h-\tilde H_\varepsilon h)(x)&\geq  
    (\varepsilon_1- \ha(1-\varepsilon_1) e^{\varepsilon} \varepsilon)\sum_{m \in \mathcal{G}_{-}} r_m \big(h(x)-h(m(x)\big)\\&+(\varepsilon_2-\ha(1+\varepsilon_2)\varepsilon) \sum_{m \in \mathcal{G}_+} r_m \big(h(m(x))-h(x)\big).
 \end{align*}
If we choose $\varepsilon_1,\varepsilon_2$ such that
 \begin{equation}\label{eqvareps}
     \frac{\varepsilon_1}{1-\varepsilon_1}=\ha e^{\varepsilon}\varepsilon  \quand \quad \frac{\varepsilon_2}{1+\varepsilon_2}=\ha \varepsilon,
 \end{equation}
then $\tilde G h-\tilde H_\varepsilon h\geq 0$ and hence $\tilde G f_{\varepsilon}\geq 0$ by the already mentioned \cite[Lemma 3]{Swa18}. By first choosing $\eps>0$ small enough and then choosing $\eps_1,\eps_2$ as in (\ref{eqvareps}), we can make  $\eps_1$ and $\eps_2$ as small as we wish, completing the proof of the claim and thereby the proof of the whole theorem.
\end{Proof}

\appendix

\section{Differentiation of the semigroup on unbounded functions}

Let $S$ be a countable set. We calculate with matrices indexed by $S$ just as we would for finite matrices. In particular, if $f\cn S\to\R$ is a function, then we write $Af(x):=\sum_{y\in S}A(x,y)f(y)$ whenever the infinite sum is well-defined, and the product of two matrices is defined in the usual way. We let $1$ denote the identity matrix, that is, $1(x,y):=1$ if $x=y$ and $:=0$ otherwise. We let $\un 1$ denote the function that is constantly one. A \emph{subprobability kernel} on $S$ is a nonnegative matrix $K$ indexed by $S$ such that $K\un 1\leq\un 1$. If this is an equality, then $K$ is a \emph{probability kernel}.

A \emph{Markov generator} on $S$ (traditionally often called \emph{Q-matrix}) is a real matrix $G$ indexed by $S$ such that $G(x,y)\geq 0$ for all $x\neq y$ and $G(x,x)=-\sum_{y:\,y\neq x}G(x,y)<\infty$ for all $x\in S$. Each Markov generator gives rise to a (possibly explosive) continuous-time Markov chain that can be constructed from its embedded discrete-time Markov chain and independent, exponentially distributed holding times. In general, the transition probabilities $P_t(x,y)$ are subprobability kernels due to the possibility of explosion. By \cite[Thms 2.8.3, 2.8.4, and 2.8.6]{Nor97}, the functions $t\mapsto P_t(x,y)$ are continuously differentiable for each $x,y\in S$, and $(P_t)_{t\geq 0}$ is given by the unique minimal solution of either of the equations
\be\label{bacfor}
{\rm(i)}\ \dif{t}P_t=GP_t\quand{\rm(ii)}\ \dif{t}P_t=P_tG\quad(t\geq 0)
\ee
with initial condition $P_0=1$. Equations (\ref{bacfor}) (i) and (ii) are Kolmogorov's \emph{backward} and \emph{forward} equation, respectively. Similarly, if  $f\cn S\to\half$ is a bounded function, then by \cite[Thm~2.7]{Swa26IPS}, $t\mapsto u_t(x):=P_tf(x)$ is continuously differentiable for all $x\in S$ and $u$ is the minimal solution of $\dif{t}u_t=Gu_t$ $(t\geq 0)$ with initial state $u_0=f$. Perhaps surprisingly, it is not in general true that $\dif{t}P_tf=P_tGf$. Indeed, one has $P_tG\un 1=0$ while $\dif{t}P_t\un 1$ is nonzero if the process is explosive.

By definition, a sequence of functions $f_n\cn S\to\R$ converges to a limit $f$ in a \emph{bounded pointwise} way, or \emph{in bp-sense} for short, if $f_n(x)\to f(x)$ for all $x\in S$ and there exists a constant $C<\infty$ such that $|f_n|\leq C$ for all $n$. We also write $\bplim_{n\to\infty}f_n=f$. Bounded pointwise convergence for functions $f_t$ depending on a real parameter $t$ is defined similarly. If $K$ is a subprobability kernel and $f_n\to f$ in bp-sense, then by dominated convergence also $Kf_n\to Kf$ in bp-sense. We will prove the following basic fact.

\bp[Differentiation of the semigroup]
Let\label{P:bdf} $(P_t)_{t\geq 0}$ be the transition kernels of a nonexplosive Markov process with countable state space $S$ and generator $G$. Let $f$ be a real function on $S$ such that $f$ and $Gf$ are bounded. Then $\bplim_{t\down 0}t^{-1}(P_t-1)f=Gf$.
\ep

The following simple lemma will be useful in what follows.

\bl[Continuous right derivative]
Assume\label{L:rdif} that $u$ and $v$ are continuous functions from $\half$ to $\R$ and that $\lim_{\eps\down 0}\eps^{-1}(u_{t+\eps}-u_t)=v_t$ $(t\geq 0)$. Then $u$ is continuously differentiable and $\dif{t}u_t=v_t$ $(t\geq 0)$.
\el

\bpro
This follows from the stronger fact, proved in \cite[Thm~IV.1.3]{Bru78}, that if the right derivative of a continuous real-valued function is continuous in a point, then the function is differentiable in that point.
\epro

Using Lemma~\ref{L:rdif} and the fact that if $K$ is a subprobability kernel and $\mbox{\rm bp-lim}_{n\to\infty}f_n=f$, then $\mbox{\rm bp-lim}_{n\to\infty}Kf_n=Kf$, it is easy to derive from Proposition~\ref{P:bdf} that
\be\label{MP}
\dif{t}P_tf=P_tGf\quad(t\geq 0)
\ee
for each real function $f$ such that both $f$ and $Gf$ are bounded. In particular, this applies when $f=\un 1$ and $Gf=G\un 1=0$, so by our earlier remark, the assumption in Proposition~\ref{P:bdf} that $G$ is nonexplosive cannot be dropped.

A \emph{Lyapunov function} for a Markov generator $G$ is a function $L\cn S\to(0,\infty)$ that satisfies $GL\leq\la L$ for some $\la\in\R$. It is not hard to see (a proof will be provided below) that this implies that
\be\label{Plat}
P_tL\leq e^{\la t}L\qquad(t\geq 0).
\ee
We cite the following result from \cite[Thm~2.12]{Swa26IPS}.

\bt[Sufficient conditions for nonexplosiveness]
Let\label{T:nonexpl} $G$ be a Markov generator on a countable set $S$. Assume that there exists a Lyapunov function $L$ for $G$ such that
\be\label{nonexpl}
\sup\big\{\!-G(x,x):L(x)\leq N\big\}<\infty\quad\forall N<\infty.
\ee
Then the continuous-time Markov chain with generator $G$ is nonexplosive.
\et

Let $(P_t)_{t\geq 0}$ be the transition kernels of a Markov process with countable state space $S$ and generator $G$. Assume that $L$ is a Lyapunov function for $G$. We say that a function $f\cn S\to\R$ is \emph{$L$-bounded} if $|f|\leq CL$ for some $C<\infty$. By (\ref{Plat}), this implies $|P_tf|\leq P_t|f|\leq Ce^{\la t}L$, so $P_tf$ is well-defined if $f$ is $L$-bounded and $P_t$ maps the space of $L$-bounded functions into itself. By definition, a sequence of functions $f_n\cn S\to\R$ converges to a limit $f$ in an \emph{$L$-bounded pointwise} way, or \emph{in $L$-bp-sense} for short, if $f_n(x)\to f(x)$ for all $x\in S$ and there exists a constant $C<\infty$ such that $|f_n|\leq C L$ for all $n$. We also write $L\mbplim_{n\to\infty}f_n=f$. By (\ref{Plat}) and dominated convergence, if $f_n\to f$ in $L$-bp-sense, then $P_tf_n\to P_tf$ in $L$-bp-sense for all $t\geq 0$. Below is the main result of this appendix. Note that we do not a priori assume that $G$ is nonexplosive and that in the special case that $L=\un 1$, the theorem follows from Proposition~\ref{P:bdf} and Theorem~\ref{T:nonexpl}.

\bt[Differentiation of unbounded functions]
Let\label{T:Lbp} $(P_t)_{t\geq 0}$ be the transition kernels of a Markov process with countable state space $S$ and generator $G$. Assume that $L$ and $L'$ are Lyapunov functions for $G$ such that
\be\label{GP}
\sup\big\{\!-G(x,x):L'(x)\leq NL(x)\big\}<\infty\quad\forall N<\infty.
\ee
Then $L\mbplim_{t\down 0}t^{-1}(P_t-1)f=Gf$ for each real function $f$ on $S$ such that $f$ and $Gf$ are $L$-bounded. 
\et

We now provide the proofs of Proposition~\ref{P:bdf}, Theorem~\ref{T:Lbp}, and formula (\ref{Plat}).\med

\bpro[of Proposition~\ref{P:bdf}]
Fix $z\in S$, let $(X_t)_{t\geq 0}$ be the Markov process with initial state $X_0=z$ and generator $G$, let $\De\sub S$ be finite, and let $\tau:=\inf\{t\geq 0:X_t\not\in\De\}$. Then the stopped process $\ti X_t:=X_{t\wedge\tau}$ $(t\geq 0)$ is the Markov process with generator $\ti Gf(x):=1_\De Gf$. Letting $\|f\|:=\sup_{x\in S}|f(x)|$ denote the supremum-norm, we have that $\|\ti Gf\|\leq 2C\|f\|$ where $C:=\sup_{x\in\De}(-G(x,x))$, so $\ti G$ is a bounded operator. It follows that its semigroup $(\ti P_t)_{t\geq 0}$ is given by
\be
\ti P_tf:=\sum_{n=0}^\infty\frac{1}{n!}t^n\ti G^nf\qquad(t\geq 0),
\ee
where the sequence is absolutely summable in the supremum-norm for each bounded $f$. Using this, it is easy to check that for each bounded real function $f$ on $S$ and $x\in S$, the function $t\mapsto\ti P_tf$ is continuously differentiable and satisfies $\dif{t}\ti P_tf(x)=\ti P_t\ti Gf(x)$ $(t\geq 0)$.

Now assume that both $f$ and $Gf$ are bounded. Let $\De_n\sub S$ be finite sets that increase to $S$ and let $\tau_n:=\inf\{t\geq 0:X_t\not\in\De_n\}$. Applying what we have just proved to the process $(X_t)_{t\geq 0}$ started in $z$ and stopped at time $\tau_n$, we see that
\be\label{stopped}
\E^z\big[f(X_{t\wedge\tau_n})\big]=f(z)+\int_0^t\!\di s\,\E^z\big[1_{\{s<\tau_n\}}Gf(X_{s\wedge\tau_n})\big]\qquad(t\geq 0,\ n\geq 0).
\ee
By assumption $(X_t)_{t\geq 0}$ is nonexplosive, so $\tau_n\to\infty$. Since both $f$ and $Gf$ are bounded, the random variable $\om\mapsto f(X_{t\wedge\tau_n}(\om))$ and the function $(\om,s)\mapsto 1_{\{\tau_n(\om)<s\}}Gf(X_{s\wedge\tau_n}(\om))$ converge to $\om\mapsto f(X_t(\om))$ and $(\om,s)\mapsto Gf(X_s(\om))$ in a bounded pointwise way, which allows us to take the limit in (\ref{stopped}) to conclude that
\be\label{MP2}
\E^z\big[f(X_t)\big]=f(z)+\int_0^t\!\di s\,\E^z\big[Gf(X_s)\big]\qquad(t\geq 0, \ z\in S).
\ee
Note that we can differentiate this formula to obtain (\ref{MP}), which we have presented as a consequence of Proposition~\ref{P:bdf}. We will show that conversely (\ref{MP}) in its integrated form (\ref{MP2}) implies Proposition~\ref{P:bdf}. Indeed, (\ref{MP2}) implies that
\be\label{int01}
t^{-1}(P_t-1)f(z)=t^{-1}\int_0^t\!\di s\,P_sGf(z)=\int_0^1\!\di s\,\E^z\big[Gf(X_{st})\big]\qquad(0<t\leq 1,\ z\in S).
\ee
Since $Gf$ is bounded, the functions $(\om,s)\mapsto Gf(X_{st}(\om))$ converge as $t\to 0$ in a bounded pointwise way to the constant $Gf(z)$. This allows us to take the limit in (\ref{int01}) to see that $t^{-1}(P_t-1)f$ converges pointwise to $Gf$. Since $Gf$ is bounded, the right-hand side of (\ref{int01}) is bounded uniformly in $t$ and $z$, so the convergence happens in a bounded pointwise way.
\epro

\bpro[of Theorem~\ref{T:Lbp} and formula (\ref{Plat})]
We will derive Theorem~\ref{T:Lbp} from Proposition~\ref{P:bdf} using a variation of Doob's h-transform. Assume that $L\cn S\to(0,\infty)$ satisfies $GL\leq\la L$ for some $\la\in\R$. Let $\ti G$ be defined by
\be\label{tiG}
\ti G(x,y):=L(x)^{-1}G(x,y)L(y)-\la1(x,y)\qquad(x,y\in S).
\ee
Then it is straightforward to check that
\be\label{GLf}
G(Lf)=L\ti Gf+\la Lf
\ee
for each $f\cn S\to\half$. It follows that $L\ti G\un 1=GL-\la L\leq 0$. This allows us to extend $\ti G$ to a Markov generator on $\ti S:=S\cup\{\dgg\}$, where $\dgg$ is a ``cemetery state'' not contained in $S$, by setting
\be
\ti G(x,\dgg):=-\sum_{y\in S}\ti G(x,y)\quad(x\in S)
\quand
\ti G(\dgg,y):=0\quad(y\in\ti S).
\ee
Let $(\ti P_t)_{t\geq 0}$ be the transition probabilities of the (possibly explosive) continuous-time Markov chain with state space $\ti S$ and generator $\ti G$. Each function $f\cn S\to\half$ can naturally be extended to a function on $\ti S$ by setting $f(\dgg):=0$. Since $\dgg$ is a trap of the process with generator $\ti G$, the semigroup $(\ti P_t)_{t\geq 0}$ maps the class of functions that are zero in $\dgg$ into itself. Identifying a function on $S$ with its natural extension to $\ti S$, we then claim that
\be\label{PLP}
P_t(Lf)=e^{\la t}L\ti P_tf\qquad(t\geq 0)
\ee
for each $f\cn S\to\half$. Note that this implies in particular that $P_tL=e^{\la t}L\ti P_t\un 1\leq e^{\la t}L$, proving formula (\ref{Plat}). It suffices to prove (\ref{PLP}) for functions $f$ that are finitely supported. For such functions both $Lf$ and $f$ are bounded functions, so the functions
\be
u_t:=P_t(Lf)\quand\ti u_t:=\ti P_tf\qquad(t\geq 0)
\ee
are the unique minimal solutions of the equations
\be\label{difu}
{\rm(i)}\ \dif{t}u_t=Gu_t\quand{\rm(ii)}\ \dif{t}\ti u_t=\ti G\ti u_t\qquad(t\geq 0)
\ee
with initial states $Lf$ and $f$, respectively. It follows from (\ref{difu})~(ii) and (\ref{GLf}) that
\be
v_t:=e^{\la t}L\ti P_tf=e^{\la t}L\ti u_t\qquad(t\geq 0)
\ee
solves
\be\label{difv}
\dif{t}v_t=e^{\la t}L\ti G\ti u_t+\la e^{\la t}L\ti u_t=e^{\la t}G(L\ti u_t)=Gv_t\qquad(t\geq 0)
\ee
with initial state $v_0=Lf$, so $u_t\leq v_t$ by the minimality of $u$, proving the inequality $\leq$ in (\ref{PLP}). To prove the opposite inequality, we observe that by replacing $f$ by $L^{-1}f$ in (\ref{GLf}) and dividing both sides of the equation by $L$, we obtain
\be\label{tiL}
\ti G(L^{-1}f)=L^{-1}Gf-\la L^{-1}f
\ee
for each $f\cn S\to\half$. Combining this with (\ref{difu})~(i), we see that
\be
w_t:=e^{-\la t}L^{-1}P_t(Lf)=e^{-\la t}L^{-1}u_t\qquad(t\geq 0)
\ee
solves
\be
\dif{t}w_t=e^{-\la t}L^{-1}Gu_t-\la e^{-\la t}L^{-1}u_t
=e^{-\la t}\ti G(L^{-1}u_t)=\ti Gw_t\quad(t\geq 0)
\ee
with initial state $w_0=f$, so $\ti u_t\leq w_t$ by the minimality of $\ti u$, proving the inequality $\geq$ in~(\ref{PLP}).

From now on, we moreover assume that there exists a function $L'\cn S\to(0,\infty)$ such that $GL'\leq\la' L'$ for some $\la'\in\R$, and (\ref{GP}) holds. We claim that $\ti G$ is nonexplosive. We apply Theorem~\ref{T:nonexpl} with $\ti L:=L'/L$ and $\ti\la:=\la'-\la$. Using (\ref{tiL}), we obtain
\be
\ti G\ti L=\ti G(L^{-1}L')=L^{-1}GL'-\la L^{-1}L'\leq\la'L^{-1}L'-\la L^{-1}L'=\ti\la\ti L.
\ee
We observe that $\ti G(\dgg,\dgg)=0$ while $\ti G(x,x)=G(x,x)-\la$ $(x\in S)$, so
\be
\sup\big\{\!-\ti G(x,x):\ti L(x)\leq N\big\}
=0\vee\big(\la+\sup\big\{\!-G(x,x):L'(x)\leq NL(x)\big\}\big)<\infty.
\ee
This shows that Theorem~\ref{T:nonexpl} is applicable\footnote{There is a subtle point here since we are identifying functions such as $\ti L$ that are defined on $S$ with functions on $S\cup\{\dgg\}$ that are zero in $\dgg$, while Lyapunov functions according to our definition should be positive everywhere. Nevertheless, Theorem~\ref{T:nonexpl} remains true even if the function $L$ there is zero in some points, as can easily be seen from the fact that in such a case we can replace $L$ by $L+1$ and $\la$ by $\la\vee 0$ to obtain a Lyapunov function that satisfies (\ref{nonexpl}).} and hence $\ti G$ is nonexplosive.

To complete the proof, it suffices to show that if $f$ is a bounded function, then one has $L\mbplim_{t\down 0}t^{-1}(P_t-1)(Lf)=G(Lf)$ provided that $G(Lf)$ is $L$-bounded. By (\ref{GLf}), the assumption that $G(Lf)$ is $L$-bounded implies that there exists a constant $C<\infty$ such that
\be
L\big|\ti Gf+\la f\big|=\big|G(Lf)\big|\leq CL.
\ee
Since $f$ is bounded, this proves that also $\ti Gf$ is bounded. This allows us to apply Proposition~\ref{P:bdf} to conclude that $\bplim_{t\down 0}t^{-1}(\ti P_t-1)f=\ti Gf$. In particular, this says that $\dif{t}\ti P_tf\big|_{t=0}=\ti Gf$, which by (\ref{PLP}) and the calculation in (\ref{difv}) implies that $\dif{t}P_t(Lf)\big|_{t=0}=G(Lf)$. Since $t^{-1}(\ti P_t-1)f$ is bounded uniformly in $t$, formula (\ref{PLP}) implies that $t^{-1}(P_t-1)(Lf)$ is $L$-bounded uniformly in $t$, completing the proof.
\epro

\paragraph{Acknowledgements.}

Both authors were supported by grant 25-16267S of the Czech Science Foundation (GA CR). In addition, the research of C.~Campailla was funded by the Austrian Science Fund (FWF) 10.55776/\allowbreak{}PAT3123425.

\end{document}